\documentclass{article}  

\usepackage{amsmath, amssymb}  
\usepackage{graphicx}
\usepackage{booktabs}
\usepackage{algorithm, algpseudocode}
\usepackage{enumitem}
\usepackage{subcaption}  
\usepackage{bm}
\usepackage{xcolor}

\DeclareMathOperator{\cl}{cl}
\DeclareMathOperator{\conv}{conv}

\newcommand{\R}{\mathbb{R}}

\newcommand{\argmin}{\mathop{\mathrm{argmin}}}

\newcommand{\argmax}{\mathop{\mathrm{argmax}}}

\newcommand{\x}{{\bm{x}}}
\newcommand{\y}{{\bm{y}}}
\newcommand{\z}{{\bm{z}}}
\newcommand{\be}{{\bm{e}}}

\newcommand{\bu}{{\bm{u}}}

\newcommand{\ba}{{\bm{a}}}
\newcommand{\bb}{{\bm{b}}}

\newcommand{\bd}{{\bm{d}}}
\newcommand{\bq}{{\bm{q}}}

\newcommand{\bxi}{\bm{\xi}}

\newcommand{\zero}{\bm{0}}

\newtheorem{proposition}{Proposition}
\newtheorem{corollary}{Corollary}
\newtheorem{definition}{Definition}
\newtheorem{remark}{Remark}
\newtheorem{example}{Example}
\newtheorem{theorem}{Theorem}

\newtheorem{problem}{Problem}
\newtheorem{proof}{Proof}
\newcommand{\qed}{\hfill \ensuremath{\square}}

\begin{document}

\title{A truncated $\varepsilon$-subdifferential method for global DC optimization}

\author{
  Adil M. Bagirov\footnote{Centre for Smart Analytics, Federation University Australia, Victoria, Australia. \texttt{a.bagirov@federation.edu.au}}
  \and
  Kaisa Joki\footnote{University of Turku, Department of Mathematics and Statistics, FI-20014 Turku, Finland. \texttt{kaisa.joki@utu.fi}}
  \and
  Marko M. M\"akel\"a\footnote{University of Turku, Department of Mathematics and Statistics, FI-20014 Turku, Finland. \texttt{ makela@utu.fi}}
  \and
  Sona Taheri\footnote{School of Science, RMIT University, Melbourne, Australia. \texttt{sona.taheri@rmit.edu.au}}
}

\date{}  

\maketitle

\begin{abstract}
We consider the difference of convex (DC) optimization problem subject to box-constraints. Utilizing $\varepsilon$-subdifferentials of DC components of the objective, we develop a new method for finding global solutions to this problem. The method combines a local search approach with a special procedure for escaping non-global solutions by identifying improved initial points for a local search. The method terminates when the solution cannot be improved further. The escaping procedure is designed using subsets of the $\varepsilon$-subdifferentials of DC components. We compute the deviation between these subsets and determine $\varepsilon$-subgradients providing this deviation. Using these specific $\varepsilon$-subgradients, we formulate a subproblem with a convex objective function. The solution to this subproblem serves as a starting point for a local search.

We study the convergence of the conceptual version of the proposed method and discuss its implementation. A large number of academic test problems demonstrate that the method requires reasonable computational effort to find higher quality solutions than other local DC optimization methods. Additionally, we apply the new method to find global solutions to DC optimization problems and compare its performance with two benchmark global optimization solvers.
\end{abstract}


\noindent{\bf MSC classes:} 90C26, 90C59, 49J52, 65K05

\section{Introduction} \label{intro}
Consider the following difference of convex (DC) optimization problem with box-constraints:
\begin{align} \label{globprobDC}
  \begin{cases}
    \text{minimize}\quad  & f(\x)\\
    \text{subject to}     & \x \in [\ba,\bb],
  \end{cases}
\end{align}
where $f(\x) = f_1(\x) - f_2(\x),~ \ba,\bb \in \R^n, ~\ba \leq \bb$ and $f_1,~f_2: \R^n\to\R$ are convex functions. Note that then $f_ 1$ and $f_2$, and thus $f$ are all continuous. The presentation $f_1-f_2$ of the DC function $f$ is called its DC decomposition whereas $f_1$ and $f_2$ are called DC components. The problem \eqref{globprobDC} has various applications, for example, in engineering, business and machine learning \cite{Bagirov2021,bagi2016}.

Using the exact penalty function approach we can include box-constraints into the first DC component by writing $f_1(\x)+\gamma \,\max\{0,a_i-x_i, x_i-b_i,~i=1,\ldots,n\}$ which is still convex. Here $\gamma > 0$ is the penalty parameter. Therefore, without loss of generality, the problem \eqref{globprobDC} can be replaced by the following unconstrained DC problem:
\begin{align} \label{probDCuncons}
  \begin{cases}
    \text{minimize}\quad  & f(\x) = f_1(\x) - f_2(\x)\\
    \text{subject to}     & \x \in \R^n.
  \end{cases}
\end{align}
In addition, we assume that $f^* = \inf \left\{f(\x):~\x \in \R^n \right\} > -\infty$.

The problem \eqref{probDCuncons} has been studied by many researchers (see, e.g., \cite{Horst1995,Horst1999,pinter1996,Tuy1998} and several methods have been developed to solve this DC problem to global optimality \cite{Horst1999,Tuy1998}. They include a Branch-and-Bound method containing three basic operations: subdivision of simplices, estimation of lower bounds and computation of upper bounds \cite{Horst1999}. Versions of the Branch-and-Bound method differ from each other in the way these operations, especially the first operation, are implemented \cite{Horst1999,Horst1996,Tuy1998}. The outer approximation method was developed in \cite{Horst1991,Horst1987}. The decomposition and parametric right-hand-side \cite{Tuy1998}, and the extended cutting angle \cite{Ferrer2015} are other methods for solving global DC programming problems. A survey of some of these methods can be found in \cite{Horst1999,Tuy1998}.

Various local search methods, aiming to find different type of stationary points critical, Clarke stationary, etc.) of DC problems, have also been developed. These methods include the difference of convex algorithm (DCA) and its variations \cite{Antao2005,Antao2012,Artacho2020,Artacho2018,AV2020}, bundle-type methods \cite{Ackooij2021,Bagirov2011,Dolgopolik2018,Gaudioso2018,Joki2017,JokBagKarMakTah2018,Oliveira2019} as well as augmented and aggregate subgradient methods \cite{Bagtahhos2021,Bagtahjokkarmak2021}. A brief survey of most of these methods can be found in \cite{Oliveira2020}.

DC optimization problems are nonconvex and may have many local minimizers. General purpose global optimization methods are time consuming and may not be efficient if the DC optimization problem has a relatively large number of variables. Although local search methods are computationally efficient, starting from some initial point, they may end up at the closest stationary point where the value of the objective function  can be significantly different from that of the global minimizer.

In many real-world applications, local minimizers with an objective value close to that of the global minimizers may still provide satisfactory solutions to the problem. One such application is the hard partitional clustering problem, where an optimal value of its objective reflects the compactness of clusters. In \cite{BagKarTah2020}, it is shown that the local minimizers of the problem with the objective value close to that of the global minimizers provide similar compact clusters to those by the global minimizer.
This observation stimulates the development of DC optimization methods which are able to find high quality solutions in a reasonable time.

In this paper, we develop a new truncated $\varepsilon$-subdifferential (TESGO) method to globally solve the DC optimization problem \eqref{probDCuncons}. The method utilizes the DC representation of the objective function. It is a combination of a local search method and a special procedure to escape from solutions which are not global minimizers. More specifically, a local search method is used to compute a stationary point (in our case a critical point) of the DC optimization problem. Then the new escaping procedure is applied to find a better starting point for the local search if the current solution is not a global one. The escaping procedure is designed using subsets of the $\varepsilon$-subdifferentials of DC components. In this procedure, utilizing $\varepsilon$-subgradients we formulate a subproblem with a convex objective function whose solutions are used as starting points for a local search. The TESGO method terminates when the solution found by a local search method cannot be improved anymore. We study the convergence of the conceptual version of the proposed method and discuss its implementation. Using results of numerical experiments we show that TESGO requires a reasonable computational effort to find higher quality solutions than the other four local methods of DC optimization. In addition, we apply TESGO to find global solutions to DC optimization problems and compare its performance with that of two benchmark global optimization solvers.

To the best of our knowledge, this is the first attempt to design a global optimization method when $\varepsilon$-subdifferentials are used with large values of $\varepsilon$ to fulfill the global DC-optimality condition obtained in \cite{Hiriart1988}.

The rest of the paper is structured as follows. Section \ref{prelim} recalls the main notations and some basic definitions from the nonsmooth analysis. The global descent direction is defined in Section \ref{globdescent}, where it is shown that such a direction can be computed using the $\varepsilon$-sudifferentials of DC components. The conceptual method is described in Section \ref{conceptual} and its heuristic version is presented in Section \ref{method}. In Section \ref{implementmethod}, the implementation of this method is discussed. Computational results and comparison of TESGO with other methods are studied in Section \ref{computational}. Finally, Section \ref{conclusions} concludes the paper.

\section{Preliminaries} \label{prelim}
In what follows, we denote by $\R^n$ the $n$-dimensional Euclidean space and the vectors of this space with boldface lowercase letters.  Furthermore, $\langle \y,\z \rangle =\sum_{i=1}^n y_i z_i$ is the inner product of vectors $\y,\z \in \R^n$, and $\|\cdot\|$ is the associated Euclidean norm. We denote by $S_1 = \{\bd\in\R^n: ~\|\bd\|= 1\}$ the unit sphere and by $B_\varepsilon(\x)$ an open ball with the radius $\varepsilon>0$ centred at $\x \in \R^n$. The convex hull of a set is denoted by ``$\conv$" and the closure of a set by ``$\cl$".

The function $f:\R^n \rightarrow \R$ is called \textit{locally Lipschitz continuous} on $\R^n$ if at any point $\x\in\R^n$ there exist a constant $L>0$ and a scalar $\varepsilon>0$ such that
$$|f(\y)-f(\z)|\leq L\|\y-\z\| \quad \text{for all}~\y,\z\in B_{\varepsilon}(\x).$$
A direction $\bd\in\R^n$ is a \textit{local descent direction} of $f:\R^n \rightarrow \R$ at a point $\x\in\R^n$ if there exists $\bar\alpha>0$ such that $f(\x+\alpha\bd)<f(\x)$ for all $\alpha\in(0,\bar\alpha]$.

The \textit{directional derivative} of $f: \R^n \rightarrow \R$ at $\x\in\R^n$ in the direction $\bd\in\R^n$ is
\begin{align*}
f'(\x;\bd)=\lim_{t\downarrow 0}\,\, \frac{f(\x+t\bd)-f(\x)}{t},
\end{align*}
if it exists. For a finite valued convex function $f:\R^n \rightarrow \R$ the directional derivative exits at any $\x\in\R^n$ in every direction $\bd\in\R^n$ and it satisfies \cite{BagKarMak2014,Rockafellar1970}
\begin{align*}
f'(\x;\bd)=\inf_{t>0}\,\, \frac{f(\x+t\bd)-f(\x)}{t}.
\end{align*}
The \textit{$\varepsilon$-directional derivative} for a convex function $f: \R^n \rightarrow \R$ at a point $\x\in\R^n$ in the direction $\bd\in\R^n$ is defined as \cite{BagKarMak2014,Demyanov1982}
\begin{align} \label{edirder}
f_{\varepsilon}'(\x;\bd)=\inf_{t>0} \,\,\frac{f(\x+t\bd)-f(\x)+\varepsilon}{t}.
\end{align}

The \textit{subdifferential} of a convex function $f:\R^n \rightarrow \R$ at a point $\x \in \R^n$ is  \cite{BagKarMak2014,Rockafellar1970}
$$\partial f(\x) = \Big\{\bxi \in \R^n:~f(\y) \geq f(\x)+\langle \bxi, \y-\x \rangle,~\forall\, \y\in\R^n \Big\}.$$
Each vector $\bxi\in\partial f(\x)$ is called a \textit{subgradient} of $f$ at $\x$. The subdifferential $\partial f(\x)$ is a nonempty, convex and compact set such that $\partial f(\x)\subseteq B_L(\zero)$, where $L>0$ is the Lipschitz constant of $f$ at $\x$.

For $\varepsilon \geq 0$, the \textit{$\varepsilon$-subdifferential} of a convex function $f:\R^n\rightarrow \R$ at a point $\x \in \R^n$ is given with  \cite{Rockafellar1970}
\begin{align*}
\partial_\varepsilon f(\x) = \Big\{\bxi_\varepsilon\in \R^n:~f(\y)\geq f(\x)+\langle \bxi_\varepsilon, \y-\x \rangle -\varepsilon,~\forall\, \y\in\R^n \Big\}.
\end{align*}
Each vector $\bxi_\varepsilon \in \partial_\varepsilon f(\x)$ is called an \textit{$\varepsilon$-subgradient} of $f$ at $\x$. The set $\partial_\varepsilon f(\x)$ is nonempty, convex and compact, and it contains the subgradient information from some neighbourhood of $\x$. This is due to the fact that $\partial f(\y) \subseteq \partial_\varepsilon f(\x)$ for all $\y \in B_{\frac{\varepsilon}{2L}}(\x)$, where $L>0$ is the Lipschitz constant of $f$ at $\x$ (see Theorem 2.33 in \cite{BagKarMak2014}).

For a locally Lipschitz continuous function $f:\R^n\rightarrow\R$, the \textit{Clarke subdifferential} at a point $\x\in\R^n$ is defined as \cite{Clarke1983}
\begin{align*}
\partial_C f(\x)=\conv\left\{ \lim_{i\rightarrow\infty} \nabla f(\x_i)\,:\,\x_i\rightarrow\x \text{ and } \nabla f(\x_i) \text{ exists} \right\},
\end{align*}
and, similarly to the convex case, each $\bxi\in\partial_C f(\x)$ is called a subgradient. Moreover, $\partial_C f(\x)=\partial f(\x)$ for a convex function $f$.

The \textit{Goldstein $\varepsilon$-subdifferential} of a locally Lipschitz continuous function $f:\R^n\rightarrow\R$ at a point $\x\in\R^n$ for $\varepsilon\geq 0$ is defined by \cite{Goldstein1977}
\begin{align*}
\partial_{\varepsilon}^G f(\x)=\cl\, \conv ~\bigcup_{\y \,\in \,{B}_\varepsilon(\x)} \partial f(\y).
\end{align*}
The set $\partial_{\varepsilon}^G f(\x)$ coincides with $\partial_C f(\x)$ with the selection $\varepsilon=0$, and for any $\varepsilon\geq 0$ we have $\partial_C f(\x)\subseteq \partial_{\varepsilon}^G f(\x)$. Thus, the Goldstein $\varepsilon$-subdifferential can be used to approximate the Clarke subdifferential and, in the convex case, the subdifferential. Moreover, for a convex function $f$ we have $\partial_{\varepsilon}^G f(\x)\subseteq\partial_{2L\varepsilon} f(\x)$, where $L > 0$ is its Lipschitz constant at a point $\x\in\R^n$ (see Theorem 3.11 in \cite{BagKarMak2014}).

Note that the execution of a local search method is stopped when a stationary point $\x^*\in\R^n$ is reached. For the DC problem \eqref{probDCuncons}, the most common option is a \textit{critical point} fulfilling the condition $\partial f_1(\x^*)\cap \partial f_2(\x^*)\neq\emptyset$. The set of critical points contains all global and local minimizers as well as \textit{Clarke stationary points} satisfying the condition $\zero\in\partial_C f(\x)$ at a point $\x\in\R^n$. In addition, \textit{inf-stationary points} fulfilling the condition $\partial f_2(\x)\subseteq \partial f_1(\x)$ at $\x\in\R^n$ are included in the set of critical points. Moreover, at a local minimizer $\x^*\in\R^n$, criticality, Clarke stationarity and inf-stationarity are satisfied. In what follows, we use the term critical point whenever we talk about stationary points.

\paragraph{Quality of solutions.} Since DC optimization problems are nonconvex, they may have many local minimizers including global ones. Some of the local minimizers are closer (in the sense of the objective value) to the global minimizer than others. This raises the question about the accuracy of local minimizers with respect to global minimizers (or the best known local minimizers).

Let $U$ be a set of critical points and $U^*$ be a set of global minimizers (or best known local minimizers) of the problem \eqref{probDCuncons}. It is clear that $U^* \subseteq U$. Set $f^*=f(\x^*)$ for $\x^* \in U^*$ and take any $\x \in U$. Then the \textit{accuracy} (\textit{relative error}) of the critical point $\x$ is defined as
\begin{equation} \label{accuracy}
E(\x) = \frac{f(\x) - f^*}{|f^*|+1}.
\end{equation}
It is clear that $E(\x) \geq 0$. We say that the critical point $\x \in U$ is \textit{higher quality} than the critical point $\y \in U$ if $E(\x) < E(\y)$. For a given $\tau > 0$, the critical point $\x \in U$ is called the $\tau$-\textit{approximate global minimizer} of the problem \eqref{probDCuncons} if $E(\x) \leq \tau$.

\section{Global optimality condition and escaping procedure} \label{globdescent}
Consider the problem \eqref{probDCuncons}. The following theorem and corollary on global optimality conditions were established in \cite{Hiriart1988}.

\begin{theorem} \label{hiriart}
For $\x^* \in \R^n$ to be a global minimizer of the problem \eqref{probDCuncons} it is necessary and sufficient that
\begin{equation} \label{hiriart1}
\partial_\varepsilon f_2(\x^*) \subseteq \partial_\varepsilon f_1(\x^*)\quad\text{for all}~\varepsilon \geq 0.
\end{equation}
\end{theorem}

\begin{corollary} \label{hiriart2}
Let $\x^* \in \R^n$ be a global minimizer of the problem \eqref{probDCuncons} and $f^*=f(\x^*)$. Let also $\bar{\x} \in \R^n$ be a local minimizer of the problem \eqref{probDCuncons} and $\bar{f} = f(\bar{\x})$. Set $\alpha = \bar{f} - f^* \geq 0$. Then
$$
\partial_\varepsilon f_2(\bar{\x})  \subseteq \partial_{\varepsilon+\alpha} f_1(\bar{\x})\quad\text{for all}~\varepsilon \geq 0.
$$
\end{corollary}

From the definition of the $\varepsilon$-approximate global minimizers and Corollary \ref{hiriart2} we get the following corollary.

\begin{corollary} \label{hiriart3}
Let $\bar{\varepsilon} \geq 0$ and $\bar{\x} \in U$ be an $\bar{\varepsilon}$-approximate global minimizer. Then
$$
\partial_\varepsilon f_2(\bar{\x}) \subseteq \partial_{\varepsilon+\hat{\varepsilon}} f_1(\bar{\x})\quad\text{for all}~\varepsilon \geq 0
$$
where $\hat{\varepsilon} = \bar{\varepsilon}(|f^*|+1)$.
\end{corollary}

Theorem \ref{hiriart} implies that if the point $\bar{\x} \in \R^n$ is a local minimizer but not a global one, then the condition \eqref{hiriart1} is not satisfied for some $\varepsilon > 0$ and there exists $\y \in \R^n$ such that $f(\y) < f(\bar{\x})$. This leads to the following definition of the global descent direction.

\begin{definition} \label{descentdir}
Let $\bar{\x} \in \R^n$ be a local minimizer of the problem \eqref{probDCuncons}. A direction $\bd \in \R^n, \bd \neq \zero$ is called a \emph{global descent direction} of the function $f$ at the point $\bar{\x}$ if there exist $\bar{\alpha}>0$ and $\delta\in(0,\bar{\alpha})$ such that $f(\bar{\x}+\alpha \bd) < f(\bar{\x})$ for all $\alpha \in (\bar{\alpha}-\delta,\bar{\alpha}+ \delta)$.
\end{definition}

Note that the global descent directions are defined at local minimizers. This is due to the fact that at any other point, there always exists local descent directions. However, local descent directions do not exist at local minimizers.

Next, we prove that if at the local minimizer $\bar{\x} \in \R^n$ the optimality condition \eqref{hiriart1} is not satisfied, $\varepsilon$-subdifferentials of DC components can be used to compute global descent directions at this point.

\begin{theorem} \label{globdescent2}
Let $\bar{\x} \in \R^n$ be a local minimizer of the problem \eqref{probDCuncons}. If the condition \eqref{hiriart1} is satisfied at $\bar{\x} \in \R^n$, then it is a global minimizer. Otherwise, there exists a global descent direction at this point.
\end{theorem}
\begin{proof} If the condition \eqref{hiriart1} is satisfied at the local minimizer $\bar{\x} \in \R^n$ for any $\varepsilon \geq 0$, then according to Theorem \ref{hiriart} $\bar{\x}$ is a global minimizer.

Now assume that a point $\bar{\x} \in \R^n$ is a local minimizer but not a global one. Then there exists $\bar{\varepsilon} > 0$ such that
$$
\partial_{\bar{\varepsilon}} f_2(\bar{\x}) \nsubseteq \partial_{\bar{\varepsilon}} f_1(\bar{\x}),
$$
as a local minimizer is always an inf-stationary point. This means that there exists $\bar{\bxi}_2 \in \partial_{\bar{\varepsilon}} f_2(\bar{\x})$ such that $\bar{\bxi}_2 \notin \partial_{\bar{\varepsilon}} f_1(\bar{\x})$. Construct the convex function
\begin{equation} \label{hatfunc}
\hat{f}(\y) = f_1(\y) - \Big[f_2(\bar{\x}) + \langle \bar{\bxi}_2, \y - \bar{\x} \rangle - \bar{\varepsilon}\Big],\quad \y \in \R^n.
\end{equation}
Then we have
$$
\hat{f}(\bar{\x}) = f_1(\bar{\x}) - f_2(\bar{\x}) + \bar{\varepsilon},
$$
or
\begin{equation} \label{funchat}
f(\bar{\x}) = \hat{f}(\bar{\x}) - \bar{\varepsilon}.
\end{equation}
Note that, the function $\hat{f}$ can be represented as a sum of two convex functions: $f_1$ and $h(\y) = -f_2(\bar{\x}) - \langle \bar{\bxi}_2, \y - \bar{\x} \rangle + \bar{\varepsilon}$. Since ${h}$ is a linear function of $\y$ for any $\varepsilon \geq 0$ we have
$$
\partial_\varepsilon h(\y) = \Big\{-\bar{\bxi}_2\Big\}.
$$
Then applying the formula for the $\varepsilon$-subdifferential of the sum of two convex functions \cite{Demyanov1982} and taking into account that for a convex  function $f$ we have $\partial_{\varepsilon_1} f(\x) \subseteq \partial_ \varepsilon f(\x)$ for all $0\leq \varepsilon_1\leq \varepsilon$ (see Theorem 2.32 (ii) in \cite{BagKarMak2014}), we obtain
\begin{equation*}
\begin{array}{lrl}
\partial_\varepsilon \hat{f}(\y) & = & \conv \bigcup\limits_{\substack{\varepsilon_1+\varepsilon_2 = \varepsilon \\\varepsilon_1, \varepsilon_2 \geq 0}} \Big[\partial_{\varepsilon_1} f_1(\y) + \partial_{\varepsilon_2} {h}(\y) \Big] \\
                        & = & \conv \bigcup\limits_{\varepsilon_1 \in [0,\varepsilon]} \Big[\partial_{\varepsilon_1} f_1(\y) - \bar{\bxi}_2\Big]\\
                        & = & \partial_\varepsilon f_1(\y) - \bar{\bxi}_2.
\end{array}
\end{equation*}
This means that the $\varepsilon$-subdifferential of the function $\hat{f}$ at a point $\y \in \R^n$ is
\begin{equation}
\partial_\varepsilon \hat{f}(\y) = \conv \Big\{\bxi \in \R^n:~~\bxi_1 \in \partial_\varepsilon f_1(\y),\,\,\bxi=\bxi_1 - \bar{\bxi}_2 \Big\}.
\end{equation}
Since $\bar{\bxi}_2 \notin \partial_{\bar{\varepsilon}} f_1(\bar{\x})$ it follows that $\zero \notin \partial_{\bar{\varepsilon}} \hat{f}(\bar{\x})$.

In addition, from $\bar{\bxi}_2 \in \partial_{\bar{\varepsilon}} f_2(\bar{\x})$ we get
$$
f_2(\y) \geq f_2(\bar{\x}) + \langle \bar{\bxi}_2, \y - \bar{\x} \rangle - \bar{\varepsilon}\quad \text{for all}~\y\in\R^n,
$$
and therefore,
\begin{equation} \label{overest1}
f(\y) = f_1(\y) - f_2(\y) \leq f_1(\y) - \Big[f_2(\bar{\x}) + \langle \bar{\bxi}_2, \y - \bar{\x} \rangle - \bar{\varepsilon}\Big] = \hat{f}(\y)\quad\text{for all}~\y \in \R^n.
\end{equation}
This together with \eqref{funchat} implies that
\begin{equation} \label{overest}
f(\y) - f(\bar{\x}) \leq \hat{f}(\y) - \hat{f}(\bar{\x}) + \bar{\varepsilon}\quad \text{for all }\y \in \R^n.
\end{equation}
Next, we compute
$$
\bar{\bxi} = \argmin_{\bxi \in \partial_{\bar{\varepsilon}} \hat{f}(\bar{\x})} ~\|\bxi\|.
$$
Since $\zero \notin \partial_{\bar{\varepsilon}} \hat{f}(\bar{\x})$ we have $\|\bar{\bxi}\| > 0$. Applying the necessary and sufficient condition for optimality (Lemma 5.2.6 in \cite{MakNei1992}) for any $\bxi \in \partial_{\bar{\varepsilon}} \hat{f}(\bar{\x})$ we obtain
$$
\langle \bar{\bxi}, \bxi - \bar{\bxi} \rangle \geq 0,
$$
which can be rewritten as
$$
\langle -\bar{\bxi}, \bxi \rangle \leq -\|\bar{\bxi}\|^2.
$$
The $\varepsilon$-directional derivative of the function $\hat{f}$ at the point $\bar{\x}$ in the direction $\bd$ with the selection $\varepsilon = \bar{\varepsilon}$ is (Theorem 2.32 in \cite{BagKarMak2014})
$$
\hat{f}^\prime_{\bar{\varepsilon}}(\bar{\x},\bd) = \max_{\bxi \in \partial_{\bar{\varepsilon}} \hat{f}(\bar{\x})} \langle \bxi, \bd \rangle.
$$
Then for $\bd = - \bar{\bxi}$ we have
\begin{equation} \label{hatder}
\hat{f}^\prime_{\bar{\varepsilon}} (\bar{\x},\bd) \leq -\|\bar{\bxi}\|^2.
\end{equation}
On the other hand, it follows from \eqref{edirder} that
\begin{equation*}
\hat{f}^\prime_{\bar{\varepsilon}}(\bar{\x},\bd) = \inf_{\alpha > 0} \frac{\hat{f}(\bar{\x}+\alpha \bd) - \hat{f}(\bar{\x}) + \bar{\varepsilon}}{\alpha}.
\end{equation*}
This together with \eqref{hatder} implies that there exists $\bar{\alpha} > 0$ such that
\begin{equation}\label{estimateS}
\hat{f}(\bar{\x}+\bar{\alpha} \bd) - \hat{f}(\bar{\x}) + \bar{\varepsilon} < -\bar{\alpha} \|\bar{\bxi}\|^2/2.
\end{equation}
Then applying \eqref{overest} for $\y = \bar{\x}+\bar{\alpha} \bd$ we obtain
$$
f(\bar{\x}+\bar{\alpha} \bd) - f(\bar{\x}) < -\bar{\alpha} \|\bar{\bxi}\|^2/2 < 0.
$$
Since the function $f$ is continuous there exists $\delta \in(0,\bar{\alpha})$ such that
$$
f(\bar{\x}+\alpha \bd) - f(\bar{\x}) < -\alpha \|\bar{\bxi}\|^2/4 < 0
$$
for all $\alpha \in (\bar{\alpha} -\delta, \bar{\alpha} + \delta)$. This means that the direction $\bd = - \bar{\bxi}$ is the direction of global descent.  \qed \end{proof}

Theorem \ref{globdescent2} shows that if at a local minimizer $\bar{\x} \in \R^n$ (not a global one) for some $\bar{\varepsilon}> 0$ there exists $\bar{\bxi}_2 \in \partial_{\bar{\varepsilon}} f_2(\bar{\x})$ such that $\bar{\bxi}_2 \notin \partial_{\bar{\varepsilon}} f_1(\bar{\x})$, then the $\bar{\varepsilon}$-subgradient  $\bar{\bxi}_2$ can be used to find a global descent direction at $\bar{\x}$. In general, there are infinitely many such $\bar{\varepsilon}$-subgradients. For example, we can define $\bar{\bxi}_2$ as follows:
$$
\bar{\bxi}_2 = \argmax_{\bxi_2 \in \partial_{\bar{\varepsilon}} f_2(\bar{\x})} \text{ min}\, \Big\{\| \bxi_1-\bxi_2\|\, :\,\bxi_1\in \partial_{\bar{\varepsilon}} f_1(\bar{\x}) \Big\}.
$$

\begin{theorem} \label{overestimate1}
Let
\begin{itemize}
\item[1.] the point $\bar{\x}$ be a local minimizer, but not a global minimizer of the function $f$;
\item[2.] $\bar{\bxi}_2\in\partial_{\bar{\varepsilon}} f_2(\bar{\x})$ and $\bar{\bxi}_2 \notin \partial_{\bar{\varepsilon}} f_1(\bar{\x})$ for some $\bar{\varepsilon} > 0$;
\item[3.] the function $\hat{f}$ be defined by \eqref{hatfunc};
\item[4.] $\bar{\bxi}_1 = \argmin\limits_{\bxi_1 \in \partial_{\bar{\varepsilon}} f_1(\bar{\x})} \big\|\bxi_1-\bar{\bxi}_2\big\|. $
\end{itemize}
Then the function $\hat{f}$ is an overestimate of the function $f$, that is $f(\y) \leq \hat{f}(\y)$ for all $\y \in \R^n$. Furthermore, the point $\bar{\x}$ is not the global minimizer of $\hat{f}$ over $\R^n$ and there exists $\bar{\alpha} > 0$ such that
\begin{equation} \label{estimate1}
\hat{f}\Big(\bar{\x}-\bar{\alpha}\big(\bar{\bxi}_1-\bar{\bxi}_2\big)\Big) - \hat{f}\big(\bar{\x}\big) \leq -\bar{\varepsilon} - \frac{\bar{\alpha}}{2} \big\|\bar{\bxi}_1 - \bar{\bxi}_2\big\|^2,
\end{equation}
and
\begin{equation} \label{estimate1S}
f\Big(\bar{\x}-\bar{\alpha}\big(\bar{\bxi}_1-\bar{\bxi}_2\big)\Big) - f\big(\bar{\x}\big) \leq  - \frac{\bar{\alpha}}{2} \big\|\bar{\bxi}_1- \bar{\bxi}_2\big\|^2.
\end{equation}
\end{theorem}
\begin{proof} The fact that the function $\hat{f}$ is an overestimate of the function $f$ follows from \eqref{overest1}. The subdifferential of the function $\hat{f}$ at a point $\y \in \R^n$ is $\partial \hat{f}(\y) =\partial f_1(\y) - \bar{\bxi}_2$ and at the point $\bar\x$ is $\partial \hat{f}(\bar{\x}) =\partial f_1(\bar{\x}) - \bar{\bxi}_2.$ Since $\bar{\bxi}_2 \notin \partial_{\bar{\varepsilon}} f_1(\bar{\x})$ and $\partial f_1(\bar{\x}) \subseteq \partial_{\bar{\varepsilon}} f_1(\bar{\x})$ it follows that $\bar{\bxi}_2 \notin \partial f_1(\bar{\x})$. This implies that $\zero \notin \partial \hat{f}(\bar{\x})$. As $\hat{f}$ is convex function we get that $\bar{\x}$ is not a global minimizer of $\hat{f}$.  The inequality \eqref{estimate1} then follows from \eqref{estimateS}. According to \eqref{overest1}, we have
$$
f\left(\bar{\x}-\bar{\alpha}(\bar{\bxi}_1 - \bar{\bxi}_2)\right) \leq \hat{f}\left(\bar{\x}-\bar{\alpha}(\bar{\bxi_1} - \bar{\bxi}_2)\right).
$$
Then applying \eqref{funchat}, the inequality \eqref{estimate1S} follows from \eqref{estimate1}. \qed \end{proof}

Theorem \ref{overestimate1} implies that if a local minimizer $\bar{\x} \in \R^n$ is not a global minimizer of the DC function $f$, then we can find $\bar{\varepsilon} > 0$ and
the $\bar{\varepsilon}$-subgradient $\bar{\bxi}_2 \in \partial_{\bar{\varepsilon}} f_2(\x)$ such that $\bar{\bxi}_2 \notin \partial_{\bar{\varepsilon}} f_1(\x)$. Then we construct the function $\hat{f}$ using the $\bar{\varepsilon}$-subgradient $\bar{\bxi}_2$ and  compute a point where the value of the function $f$ is significantly less than that of at $\bar{\x}$ by solving the problem
\begin{equation} \label{subprob3}
\mbox{minimize}~~\hat{f}(\x) ~~\mbox{subject~to}~~\x \in \R^n.
\end{equation}
This means that by solving the problem \eqref{subprob3} one can escape from the local minimizers $\bar{\x}$ and find a better starting point for a local search.

Next, using a DC function $f$ of one variable we illustrate graphs of the function $f$ and its corresponding function $\hat{f}$.

\begin{figure}[h!]
\centering
    \includegraphics[scale=.4]{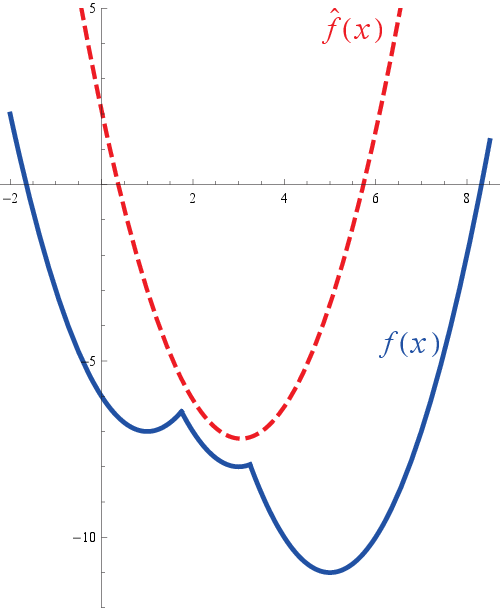}
    \caption{The graph of the function $f$ (blue) with its overestimation $\hat{f}$ (dashed red).}
    \label{functionhat}
\end{figure}

\begin{example} \label{basin7}
Consider the DC function $f(x) = f_1(x) - f_2(x), ~x \in \R,
$ where
$$
f_1(x) = x^2-5x+2,~~~f_2(x) = \max\{-3x+8, x+1, 5x-12\}.
$$
Then
$$
f(x) = \min \{x^2-2x-6, x^2-6x+1, x^2-10x+14\}.
$$
The point $\bar{x}=1$ is a local minimizer, but not a global one, of the function $f$. It can be shown that $1.1 \in \partial_{\bar{\varepsilon}} f_2(1)=[-3,13/9]$ but $1.1 \notin \partial_{\bar{\varepsilon}} f_1(1)=[-7,1]$ for $\bar{\varepsilon} = 4$. Then
$$
\hat{f}(x) = x^2 - 6.1x +2.1.
$$
It is easy to check that $f(1) = \hat{f}(1) - 4$. The graphs of the functions $f$ (blue graph) and $\hat{f}$ (dashed red graph) are depicted in Figure \ref{functionhat}. We can see  that although the function $\hat{f}$ is constructed at the point $x=1$, its minimizer is located in the more deeper basin of the function $f$.
\end{example}

\section{Conceptual algorithm} \label{conceptual}
Using results from Theorems \ref{globdescent2} and \ref{overestimate1}, we next introduce Algorithm \ref{conceptgoalg} presenting the conceptual method for globally solving the problem \eqref{probDCuncons}.

\begin{algorithm}[h]
  \caption{Conceptual algorithm} \label{conceptgoalg}
  \begin{algorithmic}[1]
   \State \label{step1} Select any starting point $\x_0 \in \R^n$ and set $k=0$.
   \State \label{step2} Apply a local search method starting from the point $\x_k$ and find a critical point $\bar{\x}_k$.
   \State \label{glopopt} If $\partial_\varepsilon f_2(\bar{\x}_k) \subseteq \partial_\varepsilon f_1(\bar{\x}_k)$ for any $\varepsilon \geq 0$, then \textbf{STOP}. The point $\bar{\x}_k$ is a global minimizer.
   \State \label{devia} Find $\bar\varepsilon \geq 0$, $\bar{\bxi}_1 \in \partial_{\bar\varepsilon} f_1(\bar{\x}_k)$ and $\bar{\bxi}_2 \in \partial_{\bar\varepsilon} f_2(\bar{\x}_k)$ such that
      $$
       \|\bar{\bxi}_1 - \bar{\bxi}_2\|^2 = \max_{\bxi_2 \in \partial_{\bar\varepsilon} f_2(\bar{\x}_k)} \min_{\bxi_1 \in \partial_{\bar\varepsilon} f_1(\bar{\x}_k)} \|\bxi_1 - \bxi_2\|^2 > 0.
      $$
   \State \label{subprobalg} Construct the function $\hat{f}(\y) = f_1(\y) - \Big[f_2(\bar{\x}_k) + \langle \bar{\bxi}_2, \y - \bar{\x}_k \rangle - \bar\varepsilon \Big]$ and solve the problem
   $$
    \mbox{minimize}~~\hat{f}(\y)~~~\mbox{subject~to}~~\y \in \R^n.
   $$
Let $\bar{\y}$ be a solution to this problem.
   \State Set $\x_{k+1} = \bar{\y}$, $k=k+1$ and go to Step \ref{step2}.
  \end{algorithmic}
\end{algorithm}

In Algorithm \ref{conceptgoalg}, Steps \ref{glopopt}$-$\ref{subprobalg} are, in general, time consuming. In Step \ref{glopopt}, we verify that the obtained critical point $\bar \x_k  \in \R^n$ satisfies global optimality. It involves the calculation of $\varepsilon$-subdifferentials of DC components for any $\varepsilon \geq 0$. This step requires the calculation of the whole $\varepsilon$-subdifferential and in general, may not be carried out in a reasonable time. In Step \ref{devia}, for a given $\bar\varepsilon\geq 0$ we compute the deviation of the set $\partial_{\bar\varepsilon} f_2(\bar \x_k)$ from the set $\partial_{\bar\varepsilon} f_1(\bar \x_k)$. If both sets are polytopes then this problem can be replaced by the finite number of quadratic programming problems which can be solved efficiently applying existing algorithms \cite{Frangioni1996,Kiwiel1989,Nurminski2008,Wolfe1976}. Finally, in Step \ref{subprobalg} we solve a subproblem to improve the current local minimizer. This problem is a convex optimization problem which can be solved efficiently.

The next theorem presents convergence result for Algorithm \ref{conceptgoalg}.

\begin{theorem} \label{convergence}
Assume that the objective function $f$ in the problem \eqref{probDCuncons} has a finite number of critical points with distinct objective values. Then Algorithm \ref{conceptgoalg} finds the global minimizer of the  problem \eqref{probDCuncons} in a finite number of iterations.
\end{theorem}

\begin{proof} Algorithm \ref{conceptgoalg} finds a new critical point at each iteration. According to \eqref{estimate1S} in Theorem \ref{overestimate1}, the value of the objective function at this critical point is strictly less than its value at the critical point found at the previous iteration. Since the number of such critical points is finite the algorithm will find the global minimizer after a finite number of iterations. \qed \end{proof}

\section{Truncated $\varepsilon$-subdifferential method} \label{method}
Step \ref{glopopt} is the most time consuming step in Algorithm \ref{conceptgoalg}. This step cannot be implemented in its present formulation as it requires the calculation of the $\varepsilon$-subdifferentials of the DC components for any $\varepsilon \geq 0.$ In order to make this step implementable, we will replace $\varepsilon$-subdifferentials by their approximations. First, we prove the following useful proposition.
\begin{proposition} \label{subset1}
Let $A,Q \subset \R^n$ be compact convex sets and for some $\delta > 0$
\begin{equation} \label{approx1}
\max_{\ba \in A} ~\langle \ba, \bd \rangle \leq \max_{\bq \in Q} ~\langle \bq, \bd \rangle + \delta,
\end{equation}
for any $\bd \in S_1$. Then
$$
A \subseteq Q + \cl\,{B}_\delta(\zero).
$$
\end{proposition}
\begin{proof} Assume on the contrary that \eqref{approx1} is true for some $\delta > 0$ but $A \not\subseteq Q + \cl\,{B}_\delta(\zero).$ This means that there exists $\bar{\ba} \in A$ such that $\bar{\ba} \not\in Q + \cl\,{B}_\delta(\zero)$ which in its turn implies that $\|\bq - \bar{\ba}\| > \delta$ for all $\bq \in Q$. Since the set $Q$ is convex and compact there exists a unique $\bar{\bq}$ minimizing the distance from $\bar{\ba}$ to the set $Q$ (see e.g. Lemma 2.2 in  \cite{BagKarMak2014}). In other words
$$
\|\bar{\bq} - \bar{\ba}\| = \min_{\bq \in Q} \|\bq - \bar{\ba}\|.
$$
It is clear that $\|\bar{\bq} - \bar{\ba}\| > \delta$. This means that $\zero \notin Q - \bar{\ba} $. Since the set $Q - \bar{\ba}$ is convex the necessary and sufficient condition (see e.g. Lemma 2.2 in \cite{BagKarMak2014}) for a minimum implies that
$$
\Big\langle \frac{ \bar{\bq} - \bar{\ba}}{\|\bar{\bq} - \bar{\ba}\|}, (\bq - \bar{\ba}) - (\bar{\bq} - \bar{\ba}) \Big\rangle \geq 0,
$$
or
\begin{equation} \label{necessary1}
\frac{\langle \bar{\bq} - \bar{\ba}, \bq - \bar{\ba} \rangle}{\|\bar{\bq} - \bar{\ba}\|} \geq \|\bar{\bq} - \bar{\ba}\|\quad\text{for all} ~\bq \in Q.
\end{equation}
Consider the following direction
$$
\bar{\bd} = -\frac{\bar{\bq} - \bar{\ba}}{\|\bar{\bq} - \bar{\ba}\|} \in S_1.
$$
Then it follows from \eqref{necessary1} that
$$
\langle \bar{\ba}, \bar{\bd} \rangle \geq \langle \bq, \bar{\bd} \rangle + \|\bar{\bq} - \bar{\ba}\| > \langle \bq, \bar{\bd} \rangle + \delta \quad\text{for all}~\bq \in Q,
$$
or
$$
\langle \bar{\ba}, \bar{\bd} \rangle > \max_{\bq \in Q} \langle \bq, \bar{\bd} \rangle + \delta.
$$
This contradicts the condition \eqref{approx1} of the proposition. \qed \end{proof}

Next, we define two sets which are used to obtain an inner and outer approximations of the $\varepsilon$-subdifferential $\partial_{\varepsilon} f(\x)$.
\begin{definition}
Let $f:\R^n\rightarrow\R$ be a convex function and $t>0$ be a given number. The set
\begin{equation} \label{setD}
D_t f(\x) = \conv \Big\{\bxi \in \R^n: ~\exists ~\bd \in S_1~~\mbox{s.t.}~~\bxi \in \partial f(\x+t\bd) \Big\}
\end{equation}
is called the \emph{$t$-spherical subdifferential} of the function $f$ at the point $\x \in \R^n$.
\end{definition}

\begin{proposition} \label{setD2}
Let $f: \R^n \to \R$ be a convex function and $t > 0$ be a given number. Then $D_t f(\x)$ is convex and compact.
\end{proposition}

\begin{proof} The convexity of the set $D_t f(\x)$ follows from its definition. For compactness, it is sufficient to show that the set
$$\bar{D}_t f(\x) = \Big\{\bxi \in \R^n: ~\exists ~\bd \in S_1~~\mbox{s.t.}~~\bxi \in \partial f(\x+t\bd) \Big\}$$ is compact. Since the subdifferential $\partial f(\x)$ is bounded on bounded sets we get that the set $\bar {D}_t f(\x)$ is bounded. Next, we show that the set $\bar{D}_t f(\x)$ is closed. Take any sequence $\{\bxi_k\},~\bxi_k \in \bar{D}_t f(\x)$. Assume that $\bxi_k \rightarrow \bar{\bxi}$ as $k \rightarrow \infty$. For each $\bxi_k$ there exists $\bd_k \in S_1$ such that $\bxi_k \in \partial f(\x+t\bd_k)$. Since the set $S_1$ is compact the sequence $\{\bd_k\}$ has at least one limit point. Without loss of generality, we assume that $\bd_k \rightarrow \bar{\bd}$. It is obvious that $\bar{\bd} \in S_1$. Upper semicontinuity of the subdifferential mapping $\x \mapsto \partial f(\x)$ (see e.g. Theorem 3.5 in \cite{BagKarMak2014}) implies that $\bar{\bxi} \in \partial f(\x+t\bar{\bd})$, and therefore $\bar{\bxi} \in \bar{D}_t f(\x)$. This means that the set $\bar{D}_t f(\x)$ is closed and together with the boundedness we get that the set is compact. Then the set $D_t f(\x)$ as the convex hull of a compact set is also compact. \qed
\end{proof}

\begin{definition}
Let $f:\R^n\rightarrow\R$ be a convex function, $L>0$ be its Lipschitz constant at a point $\x \in \R^n$, $\rho>0$ be a given number and $\delta=(2L)^{-1}\rho$. The set
\begin{equation} \label{setG}
G_\rho f(\x) = \conv ~\bigcup_{\y \,\in \,\cl\,{B}_\delta(\x)} \partial f(\y)
\end{equation}
is called a \emph{normalized Goldstein $\rho$-subdifferential} of the function $f$ at $\x \in \R^n$.
\end{definition}

\begin{remark} \label{setG1}
Let $f:\R^n\rightarrow\R$ be a convex function and $L>0$ be its Lipschitz constant at $\x\in\R^n$. Consider the function $\tilde{f}(\x) = f(\x)/L$ whose Lipschitz constant $\tilde{L}$ at $\x$ can be selected as $\tilde{L}=1$. Then the Goldstein $\rho$-subdifferential of the function $\tilde{f}$ at $\x$ contains the Goldstein $\varepsilon$-subdifferential of this function when $\rho=2\varepsilon$. Therefore, the set $G_\rho f(\x)$ is called the normalized $\rho$-subdifferential.
\end{remark}

\begin{proposition} \label{setG2}
Let $f:\R^n\rightarrow\R$ be a convex function and $\rho>0$ be a given number. Then the normalized Goldstein $\rho$-subdifferential of the function $f$ is convex and compact at any $\x \in \R^n$.
\end{proposition}

\begin{proof} The set $G_\rho f(\x)$ is convex according to its definition. Next, we show that the set
$$\bar{G}_\rho f(\x) = \bigcup_{\y \,\in \,\cl\,{B}_\delta(\x)} \partial f(\y)$$
is compact. The boundedness of this set follows from the fact that the subdifferential of $f$ is bounded on the bounded set $\cl\,{B}_\delta(\x)$.  To show the closedness of the set, take any sequence $\{\bxi_k\},~\bxi_k \in \bar{G}_\rho f(\x)$ and assume that $\bxi_k \rightarrow \bar{\bxi}$ as $k \rightarrow \infty$. For each $\bxi_k$ there exists $\y_k \in \cl B_\delta(\x)$ such that $\bxi_k \in \partial f(\y_k)$. Since the set $\cl\,B_\delta(\x)$ is compact the sequence $\{\y_k\}$ has at least one limit point. Without loss of generality, we assume that $\y_k \rightarrow \bar{\y}$. It follows from closedness of the set $\cl\,B_\delta(\x)$ that $\bar{\y} \in\cl\, B_\delta(\x)$. In addition, upper semicontinuity of the subdifferential mapping $\x \mapsto \partial f(\x)$ (see e.g. Theorem 3.5 in \cite{BagKarMak2014}) implies that $\bar{\bxi} \in \partial f(\bar{\y})$, and therefore $\bar{\bxi} \in \bar{G}_\rho f(\x)$. This means that the set $\bar{G}_\rho f(\x)$ is closed and together with boundedness compact. Then the set $G_\rho f(\x)$ as the convex hull of the compact set is also compact. \qed
\end{proof}

In the next proposition, we study the relationship between sets $D_t f(\x)$ and $G_\rho f(\x)$.

\begin{proposition}\label{l:DandG}
Let $f:\R^n \rightarrow \R$ be a convex function and $L>0$ be a Lipschitz constant of the function $f$ at a point $\x\in\R^n$. Then at $\x$ for any given $t > 0$ we have
$$
D_t f(\x) \subseteq G_\rho f(\x),
$$
when we select $\rho=2Lt$.
\end{proposition}
\begin{proof} The result follows directly from  \eqref{setD} and \eqref{setG} with the selection $\rho=2Lt$. \qed \end{proof}

Next, we show that the set $D_t f(\x)$ can be used to construct an outer approximation of the $\varepsilon$-subdifferential.
\begin{proposition} \label{p:outer}
Let $f:\R^n \rightarrow \R$ be a convex function. Then at a point $\x \in \R^n$ for any given $t>0$ and $\varepsilon \geq 0$ we have
$$
\partial_\varepsilon f(\x) \subseteq D_t f(\x) + \cl\,{B}_\delta(\zero),
$$
when we select $\delta = t^{-1} \varepsilon$.
\end{proposition}
\begin{proof} Take any $\bd \in S_1$ and any subgradient $\bar{\bxi} \in \partial f(\x + t\bd)$. Then by applying the subgradient inequality we get
$$
f(\x) - f(\x+t\bd) \geq - t\langle \bar{\bxi}, \bd \rangle,
$$
or $$f(\x+t\bd) - f(\x) \leq t\langle \bar{\bxi}, \bd \rangle.$$
Then for any $\varepsilon \geq 0$ we have
$$
f(\x+t\bd) - f(\x) + \varepsilon \leq t\langle \bar{\bxi}, \bd \rangle + \varepsilon,
$$
and thus
$$
\frac{f(\x+t\bd) - f(\x) + \varepsilon}{t} \leq \langle \bar{\bxi}, \bd \rangle + t^{-1} \varepsilon.
$$
This implies that
\begin{equation} \label{subset20S}
\begin{array}{lrl}
 f'_\varepsilon(\x,\bd) &  =   & \inf\limits_{\alpha > 0} \frac{f(\x+\alpha\bd) - f(\x) + \varepsilon}{\alpha} \\
                        & \leq & \frac{f(\x+t\bd) - f(\x) + \varepsilon}{t} \\
                        & \leq & \langle \bar{\bxi}, \bd \rangle + t^{-1} \varepsilon.
\end{array}
\end{equation}
Recall that (see Theorem 2.32 in \cite{BagKarMak2014})
$$
f'_\varepsilon(\x,\bd) = \max_{\bxi \in \partial_\varepsilon f(\x)} \langle \bxi, \bd \rangle.
$$
Then it follows from \eqref{subset20S} that for a given $\varepsilon \ge 0$
$$
\max_{\bxi \in \partial_\varepsilon f(\x)} \langle \bxi, \bd \rangle  \leq \max_{\bxi \in D_t f(\x)} \langle \bxi, \bd \rangle +t^{-1} \varepsilon.
$$
Since $\bd \in S_1$ is arbitrary, by applying Proposition \ref{subset1} we complete the proof. \qed  \end{proof}

\begin{corollary} \label{c:outer}
Let $f:\R^n \rightarrow \R$ be a convex function. Then for any given $t > 0$ at a point $\x \in \R^n$ we obtain
$$
\partial f(\x) \subseteq D_t f(\x).
$$
\end{corollary}

The following result, on the other hand, demonstrates an inner approximation of the $\varepsilon$-subdifferential.

\begin{proposition} \label{p:inner}
Let $f:\R^n \rightarrow \R$ be a convex function. Then at a point $\x\in\R^n$ for any $\rho > 0$ there exists $\varepsilon \geq 0$ such that we have
$$
G_\rho f(\x) \subseteq \partial_\varepsilon f(\x).
$$
\end{proposition}
\begin{proof} Denote $\delta=(2L)^{-1} \rho$, where $L>0$ is a Lipschitz constant of the function $f$ at $\x$. Take any $\y \in \cl\, B_\delta(\x)$ and $\bxi \in \partial f(\y)$. Then $\bxi\in G_\rho f(\x)$. The linearization error of this subgradient at $\x$ is
$$
f(\x) - f(\y) - \langle \bxi, \x - \y\rangle \geq 0.
$$
Let
$$
\varepsilon = \sup_{\substack{\y \in \cl\, B_\delta(\x)\\ \bxi \in \partial f(\y)}} ~f(\x) - f(\y) - \langle \bxi, \x - \y\rangle.
$$
Then for any $\z \in \R^n$, we have
\begin{equation} \label{subset20}
\begin{array}{lrl}
 f(\z) - f(\x) &   =  & [f(\z) - f(\y)] - [f(\x) - f(\y)] \\
               & \geq & \langle \bxi, \z - \y \rangle - [f(\x) - f(\y)] \\
               &   =  & \langle \bxi, \z - \x \rangle - [f(\x) - f(\y) - \langle \bxi, \x - \y \rangle]  \\
               & \geq & \langle \bxi, \z - \x \rangle - \varepsilon.
\end{array}
\end{equation}
This means that $\bxi \in \partial_\varepsilon f(\x)$ and the proof is  completed. \qed \end{proof}

\begin{corollary} \label{c:inner}
Let $f:\R^n \rightarrow \R$ be a convex function. Then at a point $\x\in\R^n$ for any $\varepsilon \geq 0$ we obtain
$$G_\varepsilon f(\x) \subseteq \partial_\varepsilon f(\x).$$
\end{corollary}
\begin{proof} To prove the inclusion note that (see Theorem 2.33 in  \cite{BagKarMak2014}, )
$$
\partial f(\y) \subseteq \partial_\varepsilon f(\x) \quad\text{for all}~\y \in \cl\,{B}_{\frac{\varepsilon}{2L}}(\x).
$$
Then the result follows from the definition of the set $G_\rho f(\x)$ and convexity of the $\varepsilon$-subdifferential $\partial_\varepsilon f(\x)$.  \qed \end{proof}

\begin{corollary} \label{subset7}
Let $f:\R^n \rightarrow \R$ be a convex function and $L>0$ be a Lipschitz constant of the function $f$ at a point $\x\in\R^n$. Then at $\x$ for any given $t > 0$ there exists $\varepsilon \geq 0$ such that we have
$$
D_t f(\x) \subseteq G_{2Lt} f(\x) \subseteq \partial_\varepsilon f(\x) \subseteq D_t f(\x) + \cl\,{B}_\delta(\zero),
$$
when we select $\delta=t^{-1}\varepsilon.$
\end{corollary}
\begin{proof} The first inclusion follows from Proposition \ref{l:DandG}. Other two inclusions follow from Propositions \ref{p:inner} and \ref{p:outer}, respectively. \qed \end{proof}

Corollary \ref{subset7} shows that the set $D_t f(\x)$ can be used for both inner and outer approximation of the $\varepsilon$-subdifferential $\partial_\varepsilon f(\x)$.

\begin{proposition} \label{approxcond}
Consider the problem \eqref{probDCuncons}. Assume that $D_t f_2(\x) \subseteq D_t f_1(\x)$ for any $t \geq 0$ at a point $\x\in\R^n$. Then for any $\varepsilon > 0$
$$
\partial_\varepsilon f_2(\x) \subseteq \partial_\varepsilon f_1(\x) + \cl\,{B}_\delta(\zero),
$$
where $\delta=t^{-1}\varepsilon.$
\end{proposition}
\begin{proof} Applying Corollary \ref{subset7} to the convex functions $f_1$ and $f_2$ at the point $\x$ and taking into account the condition of the proposition we have
$$
\partial_\varepsilon f_2(\x) \subseteq D_t f_2(\x)+\cl\,{B}_\delta(\zero) \subseteq D_t f_1(\x)+\cl\,{B}_\delta(\zero) \subseteq \partial_\varepsilon f_1(\x)+\cl\,{B}_\delta(\zero).
$$
This completes the proof. \qed \end{proof}

Now, we are ready to present our new truncated $\varepsilon$-subdifferential method, called TESGO, to globally solve DC problems \eqref{probDCuncons}. This method is given in Algorithm \ref{implementgoalg}. Since Proposition \ref{approxcond} implies that with some tolerance the condition ``$D_t f_2(\x) \subseteq D_t f_1(\x)$ for any $t \geq 0$" is equivalent to the global optimality condition \eqref{hiriart1}, we replace the stopping condition in Step \ref{glopopt} of Algorithm \ref{conceptgoalg} by the condition ``$D_t f_2(\x) \subseteq D_t f_1(\x)$  for any $t \geq 0$". Note also that Algorithm \ref{implementgoalg} involves the local search method and the procedure for escaping from critical points (even possibly local minimizers). The local search method is applied in Step \ref{stepin} and the escaping procedure contains Steps \ref{stopping}--\ref{subprob2}.

\begin{algorithm}[ht]
  \caption{Truncated $\varepsilon$-subdifferential (TESGO) method} \label{implementgoalg}
  \begin{algorithmic}[1]
   \State \label{input} Select a sufficiently small $\delta > 0$, any starting point $\x_0 \in \R^n$, an integer $K \geq 1$ and set $k=0$.
    \smallskip
   \State \label{stepin} Apply a local search method starting from $\x_k$ and find a critical point $\bar{\x}_k$ to solve the problem \eqref{probDCuncons}.
    \smallskip
   \State \label{boxstep} At $\bar{\x}_k$ compute
     $$
      \bar{t}_k = \max_{i=1,\ldots,n} \Big\{\bar{x}_{k,i} - a_i, b_i-\bar{x}_{k,i}\Big\},
     $$
     and set $\Delta_k = \bar{t}_k/K$ and $t_k=\Delta_k$.
   \smallskip
   \State \label{epssub} Compute the sets $D_{t_k} f_1(\bar{\x}_k)$ and $D_{t_k} f_2(\bar{\x}_k)$.
    \smallskip
   \State \label{glopopt2} If $D_{t_k} f_2(\bar{\x}_k) \subseteq D_{t_k} f_1(\bar{\x}_k) + \cl\,{B}_\delta(\zero)$, then go to Step \ref{stopping}. Otherwise, go to Step \ref{devia2}. \smallskip
   \State \label{stopping} Set $t_k=t_k + \Delta_k$. If $t_k>\bar{t}_k$, then \textbf{STOP}. The point $\bar{\x}_k$ is an approximate global minimizer. Otherwise, go to Step \ref{epssub}.
    \smallskip
   \State \label{devia2} Find $\bar{\bxi}_1 \in D_{t_k} f_1(\bar{\x}_k)$ and $\bar{\bxi}_2 \in D_{t_k} f_2(\bar{\x}_k)$ such that
      $$
       \|\bar{\bxi}_1 - \bar{\bxi}_2\|^2 = \max_{\bxi_2 \in D_{t_k} f_2(\bar{\x}_k)\,\,} \min_{\,\,\bxi_1 \in D_{t_k} f_1(\bar{\x}_k)} \|\bxi_1 - \bxi_2\|^2 > \delta.
      $$
    \smallskip
   \State \label{subprob2} Set $\varepsilon_k =  \delta t_k$, construct the function $\hat{f}_k(\y)=f_1(\y)-\Big[f_2(\bar{\x}_k) + \langle \bar{\bxi}_2, \y-\bar{\x}_k \rangle - \varepsilon_k \Big]$, solve the problem
   $$
    \mbox{minimize}~~\hat{f}_k(\y) ~~\mbox{subject~to}~~\y \in \R^n
   $$
and denote its solution by $\bar{\y}_k$.
    \smallskip
   \State \label{newit} Set $\x_{k+1} = \bar{\y}_k$, $k=k+1$ and go to Step \ref{stepin}.
  \end{algorithmic}
\end{algorithm}

Steps \ref{input}, \ref{boxstep}, \ref{stopping} and \ref{newit} of Algorithm \ref{implementgoalg} are straightforward to implement. Most time consuming steps in Algorithm \ref{implementgoalg} are Steps \ref{stepin}, \ref{epssub}, \ref{glopopt2}, \ref{devia2} and \ref{subprob2}. In Step \ref{stepin}, we apply a local search method to find a critical point of the function $f$. Here we can use any local method, and therefore this step is easily implementable. In Step \ref{epssub}, it is required to compute the sets $D_{t_k} f_1(\bar{\x}_k)$ and $D_{t_k} f_2(\bar{\x}_k)$ which is not always possible, and we discuss this in more detail in the next section. In Step \ref{glopopt2}, an approximate global optimality condition is verified. Steps \ref{glopopt2} and \ref{devia2} perform  similar tasks. Indeed, if
$$
\max_{\bxi_2 \in D_{t_k} f_2(\bar{\x}_k)\,\,\,\,} \min_{\bxi_1 \in D_{t_k} f_1(\bar{\x}_k)} \|\bxi_1 - \bxi_2\|^2 \leq \delta,
$$
then the condition in Step \ref{glopopt2} is satisfied. Both Steps \ref{glopopt2} and \ref{devia2} are implementable when the sets $D_{t_k} f_1(\bar{\x}_k)$ and $D_{t_k} f_2(\bar{\x}_k)$ are polytopes and for each vertex $\bxi_2 \in D_{t_k} f_2(\bar{\x}_k)$ we solve the quadratic programming problem
\begin{equation} \label{qpprob}
\mbox{minimize} ~~\|\bxi_1 - \bxi_2\|^2 ~~\mbox{subject~to} ~~\bxi_1 \in D_{t_k} f_1(\bar{\x}_k).
\end{equation}
There exist several algorithms developed specifically for this type of quadratic programming problem (see, for example, \cite{Frangioni1996,Kiwiel1989,Nurminski2008,Wolfe1976}). Since the number of vertices of the set $D_{t_k} f_2(\bar{\x}_k)$ is finite we have finitely many quadratic problems of the type \eqref{qpprob}, and thus,  Steps \ref{glopopt2} and \ref{devia2} can be efficiently implemented. In Step \ref{subprob2}, we solve the unconstrained convex programming problem which can be efficiently solved by existing methods.

\section{Implementation of Algorithm \ref{implementgoalg}} \label{implementmethod}
The complete calculation of the sets $D_t f_1(\x)$ and $D_t f_2(\x)$ in Step \ref{epssub} of Algorithm \ref{implementgoalg} is not always possible. However, these sets can be approximated by taking some finite point subsets of the set $S_1$. Among such subsets, positive spanning sets are the simplest and widely used ones, for example, to design direct search methods. A set of vectors $\{\bu_1,\ldots,\bu_m\},~m > 0,$ is called a positive spanning set if its positive span is $\R^n$. This set is called  positively dependent if at least one of the vectors is in the positive span generated by the remaining vectors. Otherwise, it is called positively independent. There are several well-known examples of positive spanning sets. One of them can be constructed using the standard unit vectors. Let $\{\be_1,\ldots,\be_n\}$ be the standard unit vectors in $\R^n$. Then the set
$$
U = \Big\{\pm \be_1,\ldots,\pm \be_n \Big\}
$$
is a positive spanning set in $\R^n$.

Let us take any positive spanning set
$$
U = \Big\{\bu_1,\ldots,\bu_m \Big\},\quad \text{where }\|\bu_j\| = 1,~j=1,\ldots,m,
$$
and consider the DC components $f_1$ and $f_2$ of the objective $f$. For a given $t > 0$, we construct the following sets:
\begin{equation} \label{subset4}
\tilde{D}_t f_i(\x) = \conv \Big\{\bxi \in \R^n: ~~\bxi \in \partial f_i(\x+t\bu_j),~j=1,\ldots,m \Big\},~i=1,2.
\end{equation}
It is clear that
$$
\tilde{D}_t f_i(\x) \subseteq D_t f_i(\x), ~~i=1,2.
$$
In the implementation of Step \ref{epssub} of  Algorithm \ref{implementgoalg}, we compute the sets $\tilde{D}_t f_1(\x)$ and $\tilde{D}_t f_2(\x)$ instead of the sets $D_t f_1(\x)$ and $D_t f_2(\x)$, respectively.

In numerical experiments, we consider two different versions of Algorithm \ref{implementgoalg}. The first one, called a ``simple" version, aims to significantly improve the quality of solutions obtained by a local method using limited computational effort. The second one, called a ``full" version, aims to find global solutions to DC problems using many $\varepsilon$-subgradients of DC components.

Details of the implementation of both versions of Algorithm \ref{implementgoalg} are given below:
\begin{enumerate}
\item In the problem \eqref{probDCuncons}, we fix the penalty parameter $\gamma=100$;
\item In Step \ref{input}, $\delta=0.01$ and $K = 10$ for the simple version; whereas $K=80$ for the full version of the algorithm;
\item In Step \ref{stepin}, we apply the augmented subgradient method for DC optimization (ASM-DC), introduced in \cite{Bagtahhos2021}, to find critical points. Details of the implementation of ASM-DC can be found in
\item We use the sets $\tilde{D}_t f_1(\x)$ and $\tilde{D}_t f_2(\x)$ instead of the sets $D_t f_1(\x)$ and $D_t f_2(\x)$, respectively;
\cite{Bagtahhos2021}. The maximum number of subgradient  computation $n_{max}$ at each iteration of ASM-DC is set to be $n_{max}= \max\{100, n+3\}$;
\item In Step \ref{epssub}, we compute the vertices of the sets $\tilde{D}_t f_1(\x)$ and $\tilde{D}_t f_2(\x)$ for $t=t_k$. Here, the maximum number of vertices for these sets are $m_1$ and $m_2$, respectively, where
$m_1=\min\{50,2n\}$ and $m_2=\min\{10,n\}$ for the simple version,  and  $m_1=\min\{100,2n\}$ and $m_2=\min\{30,2n\}$ for the full version.
\item In Steps \ref{glopopt2} and \ref{devia2}, for each vertex $\bxi_2 \in \tilde{D}_t f_2(\x_k)$ we apply the algorithm from \cite{Wolfe1976} to solve the quadratic programming problem \eqref{qpprob};
\item In Step \ref{subprob2}, we apply ASM-DC to solve the unconstrained convex optimization problem.
\end{enumerate}

\section{Computational results} \label{computational}
The performance of the proposed TESGO method is demonstrated by applying it to solve DC optimization test problems. Using numerical results, TESGO is compared with four local search methods of DC optimization as well as two widely used global optimization solvers. In the following, we first describe the test problems, the compared methods, and the performance measures used in our numerical experiments.  Then we discuss the results obtained.

\subsection{Test problems}
We utilize the following three groups of DC test problems:
\begin{itemize}
\item[1.] Test problems $P_1$-$P_8$ consist of the problems 2, 3, 7, 8, 10, 11, 14 and 15 described in \cite{JokBagKarMakTah2018}. They are known to have at least one or more local solutions differing from the global one;
\item[2.] Test problems $P_9$-$P_{14}$ are constructed by using the convex functions 1-3, 6-8, 10, 13, 14, 16 and 17 given in \cite{TUCreport2019}. These DC test problems are described using the notation ``Funct ~$i,j$" where $i$ and $j$ refer to the convex functions used as the first and the second DC component $f_1$ and $f_2$ of the objective, respectively. For example, the test problem $P_9$ is constructed using the convex functions 1 and 6 from \cite{TUCreport2019} as the first and second DC components, respectively;
\item[3.] Test problems $P_{15}$-$P_{20}$ are designed using some well-known global optimization test problems and modifying them as DC optimization problems. Their description is given in Appendix.
\end{itemize}

Note that in all problems from Groups 1 and 2, except $P_8$ and $P_{14}$, box-constraints are defined as  $[\ba,\bb]$, where $a_i =-100, ~b_i=100,~i=1,\ldots,n$. In $P_8$ and $P_{14}$ we have $a_i =-5, ~b_i=5,~i=1,\ldots,n$. These problems contain exponential functions, and therefore we define box-constraints for them differently to avoid very large numbers. Box-constraints for problems from Group 3 are given in their description in Appendix.

In test problems from Groups 1 and 2, only $P_5$ and $P_{11}$ have one nonsmooth DC component while the other component is smooth. In other problems from these groups, both DC components are nonsmooth. In Group 3, only $P_{15}$ has both DC components smooth. In all other problems from this group, one DC component is smooth and another one is nonsmooth. A brief description of test problems is given in Table \ref{table3G1}, where the following notations are used:
\begin{itemize}
\item $n$ - number of variables;
\item Ref - shows the label of the test problem from the referenced source;
\item $f^*$ - optimal or best known value.
\end{itemize}

\begin{table}[h!]
\caption{A brief description of test problems}
\begin{center}
\setlength\tabcolsep{6pt}
\resizebox{0.99\textwidth}{!}{
\begin{tabular}{@{}llrrrllrrrllrr@{}}
\toprule
Prob. & Ref. & $n$ & $f^*$  &&  Prob.  &  Ref.  &  $n$  &  $f^*$  &&  Prob.  &  Ref.  &  $n$  &  $f^*$  \\
\hline
\multicolumn{14}{c}{Problems from Group 1 }\\
\hline
$P_1$ & Prob 2  &  2 &   0.0000 && $P_5$  & Prob 10 & 100  &  -98.5000 && $P_7$ & Prob 14 & 200 &  0.0000 \\
$P_2$ & Prob 3  &  4 &   0.0000 && $P_5$  & Prob 10 & 200  & -198.5000 && $P_8$ & Prob 15 &   5 &  0.0000 \\
$P_3$ & Prob 7  &  2 &   0.5000 && $P_6$  & Prob 11 &   3  &  116.3333 && $P_8$ & Prob 15 &  10 &  0.0000 \\
$P_4$ & Prob 8  &  3 &   3.5000 && $P_7$  & Prob 14 &   2  &    0.0000 && $P_8$ & Prob 15 &  50 &  0.0000 \\
$P_5$ & Prob 10 &  2 &  -0.5000 && $P_7$  & Prob 14 &   5  &    0.0000 && $P_8$ & Prob 15 & 100 &  0.0000 \\
$P_5$ & Prob 10 &  5 &  -3.5000 && $P_7$  & Prob 14 &  10  &    0.0000 && $P_8$ & Prob 15 & 200 &  0.0000 \\
$P_5$ & Prob 10 & 10 &  -8.5000 && $P_7$  & Prob 14 &  50  &    0.0000 &&       &         &     &         \\
$P_5$ & Prob 10 & 50 & -48.5000 && $P_7$  & Prob 14 & 100  &    0.0000 &&       &         &     &         \\
\hline
\multicolumn{14}{c}{Problems from Group 2 }\\
\hline
$P_9$    & Funct 1,6 &  2 & -153.3333 && $P_{11}$ & Funct 3,8   & 10 &  -8.5000 && $P_{13}$  & Funct 13,17 & 10 &  -49.9443 \\
$P_9$    & Funct 1,6 &  5 & -436.6667 && $P_{11}$ & Funct 3,8   & 50 & -48.5000 && $P_{13}$  & Funct 13,17 & 50 & -273.6652 \\
$P_9$    & Funct 1,6 & 10 & -929.0909 && $P_{11}$ & Funct 3,8   &100 & -98.5000 && $P_{13}$  & Funct 13,17 &100 & -555.5672 \\
$P_9$    & Funct 1,6 & 50 &-4921.9608 && $P_{11}$ & Funct 3,8   &200 &-198.5000 && $P_{13}$  & Funct 13,17 &200 &-1116.3273 \\
$P_9$    & Funct 1,6 &100 &-9920.9901 && $P_{12}$ & Funct 13,10 &  2 &   0.0000 && $P_{14}$  & Funct 16,14 &  2 &   -1.0000 \\
$P_{10}$ & Funct 2,7 &  2 & -247.8125 && $P_{12}$ & Funct 13,10 &  5 &  -1.8541 && $P_{14}$  & Funct 16,14 &  5 &   -3.4167 \\
$P_{10}$ & Funct 2,7 &  5 & -578.4626 && $P_{12}$ & Funct 13,10 & 10 &  -4.9443 && $P_{14}$  & Funct 16,14 & 10 &  -11.2897 \\
$P_{10}$ & Funct 2,7 & 10 &-1006.8616 && $P_{12}$ & Funct 13,10 & 50 & -29.6656 && $P_{14}$  & Funct 16,14 & 50 & -126.9603 \\
$P_{10}$ & Funct 2,7 & 50 &-3564.2275 && $P_{12}$ & Funct 13,10 &100 & -60.5673 && $P_{14}$  & Funct 16,14 &100 & -320.7378 \\
$P_{10}$ & Funct 2,7 &100 &-7297.9530 && $P_{12}$ & Funct 13,10 &200 &-122.3707 && $P_{14}$  & Funct 16,14 &200 & -777.6051 \\
$P_{11}$ & Funct 3,8 &  2 &   -0.5000 && $P_{13}$ & Funct 13,17 &  2 &  -5.0000 &&           &             &    &           \\
$P_{11}$ & Funct 3,8 &  5 &   -3.5000 && $P_{13}$ & Funct 13,17 &  5 & -21.8541 &&           &             &    &           \\
\hline
\multicolumn{14}{c}{Problems from Group 3}\\
\hline
$P_{15}$ & Prob 1 &  2 & -0.3524 && $P_{17}$ & Prob 3 &   2 &  -0.8332 && $P_{19}$ & Prob 5 &   2 &-0.2500 \\
$P_{16}$ & Prob 2 &  2 &  0.0000 && $P_{18}$ & Prob 4 &   2 &  -0.3750 && $P_{20}$ & Prob 6 &   2 & 0.0000 \\
$P_{16}$ & Prob 2 &  5 &  0.0000 && $P_{18}$ & Prob 4 &   5 &  -1.3750 && $P_{20}$ & Prob 6 &   5 & 0.0000 \\
$P_{16}$ & Prob 2 & 10 &  0.0000 && $P_{18}$ & Prob 4 &  10 &  -3.0417 && $P_{20}$ & Prob 6 &  10 & 0.0000 \\
$P_{16}$ & Prob 2 & 50 &  0.0000 && $P_{18}$ & Prob 4 &  50 & -16.3750 && $P_{20}$ & Prob 6 &  50 & 0.0000 \\
$P_{16}$ & Prob 2 &100 &  0.0000 && $P_{18}$ & Prob 4 & 100 & -33.0417 && $P_{20}$ & Prob 6 & 100 & 0.0000 \\
$P_{16}$ & Prob 2 &200 &  0.0000 && $P_{18}$ & Prob 4 & 200 & -66.3750 && $P_{20}$ & Prob 6 & 200 & 0.0000 \\
\hline
\end{tabular}
}
\end{center}
\label{table3G1}
\end{table}

\subsection{Methods for comparison}
We present two different comparisons. In the first case, we aim to show that TESGO is able to escape from critical points found by a local search method and improve the quality of solution significantly using limited computational effort. In this case, we apply the simple version of TESGO for each test problem with 20 starting points randomly generated exploiting corresponding box-constraints. These same starting points are also used in other methods. We employ performance profiles to accomplish comparisons. The following local search methods of DC optimization are included in comparisons:
\begin{enumerate}
\item The aggregate subgradient method (AggSub) \cite{Bagtahjokkarmak2021};
\item The double bundle method (DBDC) \cite{JokBagKarMakTah2018};
\item The difference of convex algorithm (DCA) \cite{Antao2005} where the proximal bundle method, implemented in \cite{MakNei1992}, is used to solve convex subproblems;
\item The augmented subgradient method (ASM-DC) \cite{Bagtahhos2021}.
\end{enumerate}

In the second case, we evaluate the performance of TESGO as a global optimization solver. Thus, we use one starting point for each test problem and apply the full version of TESGO. Starting points for problems from Group 1 are given in \cite{JokBagKarMakTah2018} and for problems from Group 2 are selected as follows:
\begin{itemize}
\item In $P_9$ and $P_{10}$: $\x_0=(x_{0,1},\ldots,x_{0,n})$, where $x_{0,i} = i,~i=1,\ldots,\lfloor n/2 \rfloor$ and $x_{0,i} =-i,~i=\lfloor n/2 \rfloor + 1,\ldots,n$;
\item In $P_{11}$ and $P_{12}$: $\x_0=(x_{0,1},\ldots,x_{0,n})$, where $x_{0,i} = 0.5i,~i=1,\ldots,n$;
\item In $P_{13}$: $\x_0=(x_{0,1},\ldots,x_{0,n})$, where $x_{0,i} = 2,$ for even $i$ and $x_{0,i} =-1.5$, for odd $i=1,\ldots,n$;
\item In $P_{14}$: $\x_0=(x_{0,1},\ldots,x_{0,n})$, where $x_{0,i} = -1,$ for even $i$ and $x_{0,i} = 1$, for odd $i=1,\ldots,n$.
\end{itemize}
Starting points for problems from
Group 3 are selected as the center of the corresponding boxes.

The following two widely used global optimization solvers are applied for comparison:
\begin{itemize}
\item[1.] BARON \cite{BARON} version 23.1.5 of GAMS version 42.2.0;
\item[2.] LINDOGlobal \cite{LindoGlobal} version 14.0.5099.204 of GAMS version 42.2.0.
\end{itemize}
NEOS server \cite{czyzyketal1998,dolan2001,groppmore1997} version 6.0 is used to run these two global optimization solvers.

All the experiments (except those with BARON and LINDOGlobal) are carried out on an Intel$^{\circledR}$ Core\texttrademark\, i5-7200 CPU (2.50GHz, 2.70GHz) running under Windows 10. To compile the codes, we use gfortran, the GNU Fortran compiler. ASM-DC, DCA and TESGO are coded in Fortran 77 while DBDC and AggSub are coded in Fortran 95. We apply the implementations and default values of parameters of AggSub, ASM-DC, DBDC and DCA that are recommended in their references.


\subsection{Performance measures}
We apply performance profiles to compare the local search methods. For the number of function evaluations and the computational time (CPU time), we use the standard performance profiles introduced in \cite{Dolan2002}. To compare the accuracy of solutions obtained by these methods, we modify the standard performance profiles as described below.

The relative error $E_s(\bar{\x})$ of the solution $\bar{\x}$ obtained by the solver $s$ is defined in \eqref{accuracy}. Assume we have $k$ solvers $S=\{s_1,\ldots,s_k\}$ and the collection of $m$ problems $P=\{p_1,\ldots,p_m\}$. Applying the solver $s_i,~i=1,\ldots,k$ to the problem $p_t,~t \in \{1,\ldots,m\}$, we get solutions $\x_{1t},\ldots,\x_{kt}$ with the objective values $f(\x_{1t}),\ldots,f(\x_{kt})$. Some of these solutions may coincide. Denote by
$$
V_t = \min_{i=1,\ldots,k} f(\x_{it}),~t=1,\ldots,m.
$$
The accuracy $E_{it}$ of the solution $\x_{it}$ is defined as
$$
E_{it} = \frac{f(\x_{it}) - V_t}{|V_t|+1}.
$$
It is clear that $E_{it} \geq 0$ for any $i \in\{1,\ldots,k\}$ and $t\in\{1,\ldots,m\}$. Compute
$$
E_{max} = \max_{i=1,\ldots,k} ~\max_{t=1,\ldots,m} ~E_{it}.
$$
Take any accuracy threshold $\bar{E} \in [0,E_{max}]$. For a given solver $s_i$ and $\tau \in [0, \bar{E}]$, consider the set
$$
A_i(\tau) = \Big\{p_t ~:~E_{it} \leq \tau,~t \in \{1,\ldots,m\} \Big\},
$$
and define the following function
\begin{equation} \label{pperrors}
\sigma_i(\tau) = \frac{|A_i(\tau)|}{m}.
\end{equation}
It is clear that $\sigma_i(\tau) \in [0,1]$. The value $\sigma_i(0)$ shows the fraction of problems where the solver $s_i$ finds the best solutions. If $\sigma_i(\bar{E})=1$, the solver $s_i$ solves all problems with the given accuracy threshold.

We also use the number of function and subgradients evaluations as performance measures to compare both local and global search methods.

\subsection{Comparison with local search methods}
\begin{figure} [ht]
 \centering
   \includegraphics[width=.99\textwidth]{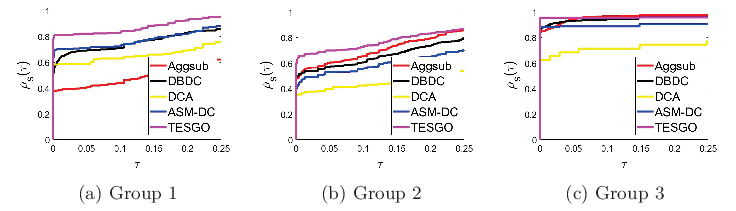}
   \caption{Performance profiles for the accuracy of solutions with $\tau=0.25$.}
\label{error123_0.25}
\end{figure}

\begin{figure} [ht]
 \centering
   \includegraphics[width=.99\textwidth]{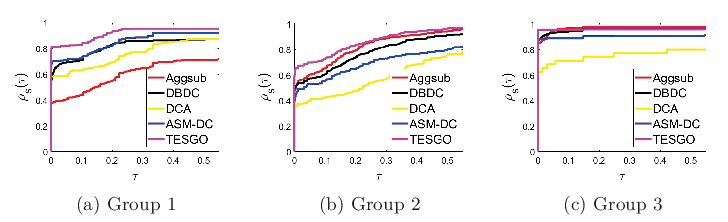}
   \caption{Performance profiles for the accuracy of solutions with $\tau=0.55$.}
  \label{error123_0.55}
\end{figure}

\begin{figure} [ht]
 \centering
   \includegraphics[width=.99\textwidth]{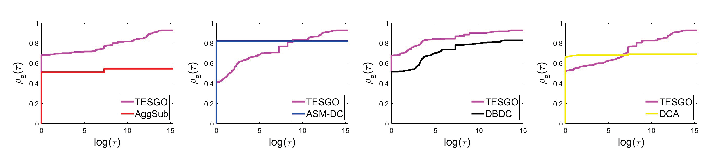}
   \caption{Performance profiles for problems from Group 1 using CPU time.}
\label{cpu1}
\end{figure}

\begin{figure} [ht]
 \centering
   \includegraphics[width=.99\textwidth]{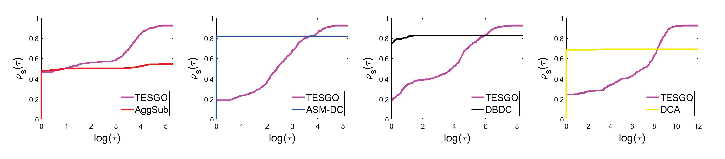}
   \caption{Performance profiles for problems from Group 1 using number of function evaluations.}
\label{nf1}
\end{figure}

\begin{figure} [ht]
 \centering
   \includegraphics[width=.99\textwidth]{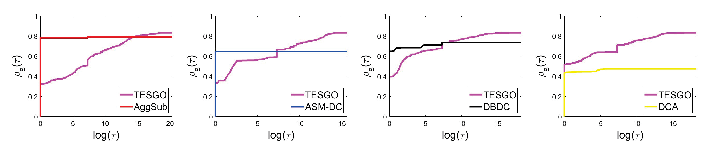}
   \caption{Performance profiles for problems from Group 2 using CPU time.}
\label{cpu2}
\end{figure}

\begin{figure} [ht]
 \centering
   \includegraphics[width=.99\textwidth]{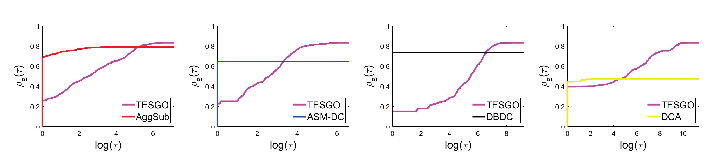}
   \caption{Performance profiles for problems from Group 2 using number of function evaluations.}
\label{nf2}
\end{figure}

\begin{figure} [ht]
 \centering
   \includegraphics[width=.99\textwidth]{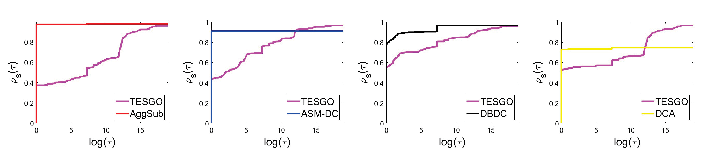}
   \caption{Performance profiles for problems from Group 3 using CPU time.}
\label{cpu3}
\end{figure}

\begin{figure} [ht]
 \centering
   \includegraphics[width=.99\textwidth]{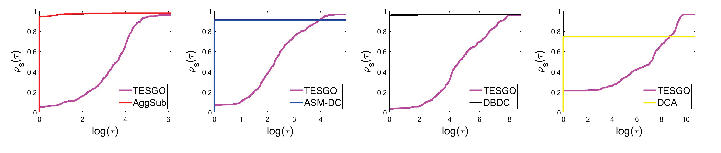}
   \caption{Performance profiles for problems from Group 3 using number of function evaluations.}
\label{nf3}
\end{figure}

We apply performance profiles using accuracy of solutions, the number of function evaluations and the CPU time to compare TESGO with local search methods. In Figures \ref{error123_0.25}--\ref{nf3}, we illustrate the results separately for each group of test problems (Group 1, Group 2 and Group 3) as described in Table \ref{table3G1}.

We say that a solver $s$ solves a problem $p$ if the solution is a $\tau$-approximate global minimizer and only such solutions are used to compute  performance profiles. Moreover, in this case $\bar{E} = E_{max}$ when the accuracy of solutions is considered. In what follows, we consider values $\tau=0.25$ and $\tau=0.55$ when performance profiles for the accuracy of solutions defined in \eqref{pperrors} are applied. These results are illustrated in Figures \ref{error123_0.25} and \ref{error123_0.55}. The higher the graph of a method, the better the method is for finding high quality solutions.

We can see from Figure \ref{error123_0.25} that in all three groups the proposed TESGO method outperforms other local methods in finding high quality solutions. Results obtained by local methods using 20 randomly generated starting points show that problems from Groups 1 and 2 have many local minimizers while problems from Group 3 have very few. In the global optimization context, this means problems from Groups 1 and 2 are more difficult than those in Group 3. In Groups 1 and 2, the difference between TESGO and local methods is considerable. In particular, in Group 1 TESGO finds the global minimizers in approximately 82\% of problems whereas the best performed local method, ASM-DC, finds such minimizers in approximately 67\% of problems. In Group 2, TESGO finds the global minimizers in almost 63\% of problems whereas the best performed local method, DBDC, finds such minimizers in 50\% of problems. Furthermore, in Groups 1 and 2, TESGO finds $\tau$-approximate global minimizers with $\tau=0.05$ in almost 83\% and 70\% of problems, respectively, whereas best performed local methods, ASM-DC and AggSub, find such solutions in almost 68\% and 58\% of problems, respectively. The TESGO method outperforms local methods also in Group 3, however the difference between TESGO and the best performed local method, DBDC is not significant. To conclude, results presented in Figures \ref{error123_0.25} and \ref{error123_0.55} show that TESGO is efficient in escaping from local minimizers and in finding high quality solutions using the limited number of $\varepsilon$-subgradients of DC components.

Next, we present performance profiles using CPU time and the number of function evaluations and provide a pairwise comparison of TESGO with AggSub, ASM-DC, DBDC and DCA. Since some of the methods use different amount of DC component evaluations, the number of function evaluations for each run of algorithms is calculated as an average of the number of evaluations of the first and second DC components.  The number of subgradient evaluations follows the similar trends with the number of the function evaluations with all the solvers, and thus, we omit these results. In addition, only results with $\tau$-approximate global minimizers with the selection $\tau=0.2$ are considered. Recall that in the standard performance profiles, the value of $\rho_s(\tau)$ at $\tau=0$ shows the ratio of the test problems for which the solver $s$ is the best --- that is, the solver $s$ uses the least computational time or evaluations --- while the value of $\rho_s(\tau)$ at the rightmost abscissa gives the ratio of the test problems that the solver $s$ can solve --- that is, the robustness of the solver $s$. In addition, the higher is a particular curve, the better is the corresponding solver. Results presented in Figures \ref{cpu1}$-$\ref{nf3} clearly show that TESGO is more robust than local methods, used in numerical experiments, across three groups of test problems. The only exception is AggSub which has a similar robustness in Group 3. On the other side, TESGO uses, in general, significantly more CPU time and function evaluations than the other local methods. This is expected as TESGO escapes from local minimizers and applies a local method multiple times.


In addition to performance profiles, we  report the number of function and subgradient evaluations of DC components in Tables \ref{NFNGDCCOMG1} --  \ref{NFNGDCCOMG3}. In these tables, $n_f$ stands for the number of evaluations of the objective function $f$ whereas $n_{f_i}$ is the number of function evaluations and $n_{g_i}$ is the number of subgradient evaluations of the DC component $f_i,~i=1,2$. We use the notation $n_{fg_i}$ when $n_{f_i}=n_{g_i},~i=1,2$. In AggSub and DBDC, we only report $n_f$ as for them we have $n_{f_1} = n_{f_2}$. In DCA, $n_{f_i} = n_{g_i},~i=1,2$ and thus we report the values of $n_{fg_1}$ and $n_{fg_2}.$  Since we use 20 starting points for each problem, we report the mean value over 20 runs of the methods.

\begin{table}[h!]
\caption{Number of function and subgradient evaluations of DC components (Group 1).}
\begin{center}
\setlength\tabcolsep{6pt}
\resizebox{0.99\textwidth}{!}{
\begin{tabular}{@{}lrrrrrlrrrrrlrrrrrlrrrrrlrrrrrlrrrrr@{}}
\toprule
Prob.       & $n$&& \multicolumn{3}{c}{AggSub}&&\multicolumn{3}{c}{DBDC} &&\multicolumn{2}{c}{DCA} &&\multicolumn{4}{c}{ASM-DC} &&\multicolumn{4}{c}{TESGO}\\
\hline
            &    && $n_f$  & $n_{g_1}$ & $n_{g_2}$ && $n_f$  & $n_{g_1}$ & $n_{g_2}$ && $n_{fg_1}$  & $n_{fg_2}$ && $n_{f_1}$  & $n_{f_2}$  & $n_{g_1}$ & $n_{g_2}$ && $n_{f_1}$  & $n_{f_2}$  & $n_{g_1}$ & $n_{g_2}$    \\
 \cline{4-6} \cline{8-10} \cline{12-13} \cline{15-18} \cline{20-23}
$P_1$	&	2	&&	199	&	101	&	40	&&	94	&	93	&	92	&&	37	&	2	&&	321	&	211	&	110	&	46	&&350  &  229   &  174	&	77	\\
$P_2$	&	4	&&4895	&4754	&	89	&&	35	&	36	&	32	&&	30	&	2	&&	630	&	391	&	240	&	71	&&4240 & 1932	& 1924	&	400	\\
$P_3$	&	2	&&3748	&3634	&	75	&&	582	&	581	&	526	&&	811	&	56	&&	420	&	272	&	148	&	49	&&4917 & 1385	& 1984	&	268	\\
$P_4$	&	3	&&	173	&	121	&	32	&&	139	&	139	&	62	&&	40	&	3	&&	385	&	257	&	128	&	42	&&3065 &  874	& 1291	&	180	\\
$P_5$	&	2	&&	127	&	82	&	27	&&	59	&	59	&	35	&&	6	&	2	&&	235	&	155	&	81	&	30	&&1733 &  606	&  704	&	140	\\
$P_5$	&	5	&&	144	&	94	&	31	&&	107	&	108	&	56	&&	7	&	2	&&	282	&	180	&	103	&	33	&&3050 & 1169	& 1339	&	282	\\
$P_5$	&	10	&&	144	&	90	&	30	&&	101	&	102	&	54	&&	8	&	2	&&	334	&	209	&	126	&	35	&&3964 & 1438	&1893	&	385	\\
$P_5$	&	50	&&	147	&	89	&	31	&&	93	&	94	&	52	&&	8	&	2	&&	541	&	317	&	224	&	38	&&3272 & 1113	&2022	&	268	\\
$P_5$	&	100	&&	146	&	98	&	33	&&	81	&	82	&	50	&&	9	&	3	&&	624	&	356	&	268	&	35	&&3531 & 1214	&2193	&	252	\\
$P_5$	&	200	&&	165	&	119	&	36	&&	60	&	61	&	40	&&	9	&	3	&&	597	&	345	&	252	&	37	&&2894 &  869	&1959	&	224	\\
$P_6$	&	3	&&4549	&4446	&	64	&&	29	&	29	&	24	&&	29	&	2	&&	348	&	220	&	128	&	41	&& 603 &  329	&309	&	100	\\
$P_7$	&	2	&&	53	&	30	&	20	&&	7	&	8	&	8	&&	16	&	2	&&	135	&	89	&	46	&	23	&& 892 &  329   &371	&	114	\\
$P_7$	&	5	&&	267	&	194	&	33	&&7864	&7863	&7862	&&	45	&	2	&&	590	&	369	&	221	&	56	&& 7722&2239	&3169	&	395	\\
$P_7$	&	10	&&	387	&	281	&	43	&&9529	&9527	&9524	&&	49	&	2	&&1220	&	757	&	463	&	88	&&13938&4466	&5697	&	635	\\
$P_7$	&	50	&&	767	&	677	&	50	&&10014	&10012	&10010	&&	59	&	2	&&3159	&1775	&1383   &	144	&&22435&5057	&10321	&	533	\\
$P_7$	&	100	&& 1231	&  1163	&	52	&&10016	&10012	&10010	&&	60	&	2	&&6308	&3403	&2905   &	201	&&29482&6198	&13918	&	495	\\
$P_7$	&	200	&&10868	& 10080	&	626	&&10022	&10015	&10014	&&	77	&	2	&&10851	&5740	&5111   &	281	&&35181&7007	&16772	&	476	\\
$P_8$	&	5	&&	831	&	792	&	45	&&1176	&1133	&   92	&&9085	&	3	&&691	&425	&	267	&	49	&&16603&3796	&7078	&	489	\\
$P_8$	&	10	&& 3003	&  2834	&	154	&&2400	&2329	&238	&&6245	&	3	&&1501	&848	&	654	&	69  &&32038&8516	&14396	&	776	\\
$P_8$	&	50	&&18423	& 18217	&	202	&& 7563	&7492	&1988	&&109471&	9	&&16343	&8369	&7975   &	202	&&76821&34120	&37888	&	909	\\
$P_8$	&  100	&&27456	&27273	&	194	&&30833	&30545	&21325	&&70195	&	9	&&45771	&23142	&22630  &	305	&&85168&37665	&42472	&	625	\\
$P_8$	&  200	&&75375	&75178	&	215	&& 5850	&5791	&5082	&&404000&	20	&&8918	&4581	&4337   &	145	&&32898&12431	&16451	&	499	\\
\hline
\end{tabular}
}
\end{center}
\label{NFNGDCCOMG1}
\end{table}

\begin{table}[h!]
\caption{Number of function and subgradient evaluations of DC components (Group 2).}
\begin{center}
\setlength\tabcolsep{6pt}
\resizebox{0.99\textwidth}{!}{
\begin{tabular}{@{}lrrrrrlrrrrrlrrrrrlrrrrrlrrrrrlrrrrr@{}}
\toprule
Prob.       & $n$&& \multicolumn{3}{c}{AggSub}&&\multicolumn{3}{c}{DBDC} &&\multicolumn{2}{c}{DCA} &&\multicolumn{4}{c}{ASM-DC} &&\multicolumn{4}{c}{TESGO}\\
\hline
            &    && $n_f$  & $n_{g_1}$ & $n_{g_2}$ && $n_f$  & $n_{g_1}$ & $n_{g_2}$ && $n_{fg_1}$  & $n_{fg_2}$ && $n_{f_1}$  & $n_{f_2}$  & $n_{g_1}$ & $n_{g_2}$ && $n_{f_1}$  & $n_{f_2}$  & $n_{g_1}$ & $n_{g_2}$    \\
 \cline{4-6} \cline{8-10} \cline{12-13} \cline{15-18} \cline{20-23}
$P_9$    &	2	&&	409	&	277	&	55	&&	42	&	35	&	29	&&	39	&	3	&&	355	&	236	&	118	&	43	&&	4588  &	 797 &  1962  &  167 \\
$P_9$    &	5	&&1141	&	870	&	113	&&	76	&	63	&	45	&&	116	&	3	&&	817	&	524	&	292	&	73	&& 10077  &	2377 &	4136  &  393 \\
$P_9$    &	10	&&2365	&1880	&	188	&&	121	&	110	&	76	&&	180	&	2	&&1392	&	851	&	541	&	97	&& 13772  &	3406 &	6030  &  488 \\
$P_9$    &	50	&&10666	&9817	&	304	&&	644	&	624	&	200	&&13606	&	2	&&6838	&3682	&3156   &	161	&& 22334  &	9858 & 11277  &  467 \\
$P_9$    &	100	&&20268	&19591	&	268	&&1589	&1561	&	559	&&31362	&	2	&&14226	&7473	&6754   &	209	&& 31013  &14248 & 15613  &  450 \\
$P_{10}$ &	2	&&357	&	229	&	53	&&	40	&	29	&	24	&&	27	&	3	&&	350	&	239	&	111	&	43	&&  4901  &	 823 &	2050  &  174 \\
$P_{10}$ &	5	&&1071	&	948	&	84	&&	62	&	55	&	44	&&	110	&	3	&&	825	&	543	&	282	&	76	&& 10594  & 2704 &	4344  &  431 \\
$P_{10}$ &	10	&&2374	&2129	&	172	&&	122	&	115	&	72	&&	257	&	3	&&1604	&	974	&	630	& 114   && 17969  & 5025 &	7817  &  665 \\
$P_{10}$ &	50	&&25614	&24782	&	721	&&1178	&1158	&	447	&&4415	&	3	&&23120	&12029	&11091	& 506	&& 61393  &30073 & 30080  & 1186 \\
$P_{10}$ &	100	&&22809	&22398	&	376	&&2865	&2844	&1121	&&13654	&	4	&&83252	&42449	&40803	&1063   && 97002  &47538 &	48209	&1282 \\
$P_{11}$ &	2	&&	127	&	82	&	27	&&	59	&	59	&	35	&&	6	&	2	&&	235	&	155	&	81	&	30	&&  4548  &595	&	2041	&	142	\\
$P_{11}$ &	5	&&	144	&	94	&	31	&&	107	&	108	&	56	&&	7	&	2	&&	282	&	180	&	103	&	33	&&  7111  &1204	&	3228	&	297	\\
$P_{11}$ &	10	&&	144	&	90	&	30	&&	101	&	102	&	54	&&	8	&	2	&&	334	&	209	&	126	&	35	&&  8598  &1480	&	4059	&	410	\\
$P_{11}$ &	50	&&	147	&	89	&	31	&&	93	&	94	&	52	&&	8	&	2	&&	541	&	317	&	224	&	38	&&  7778  &1248	&	4165	&	295	\\
$P_{11}$ &	100	&&	146	&	98	&	33	&&	81	&	82	&	50	&&	9	&	3	&&	624	&	356	&	268	&	35	&&  7500  &1181	&	4083	&	266	\\
$P_{11}$ &	200	&&	165	&	119	&	36	&&	60	&	61	&	40	&&	9	&	3	&&	597	&	345	&	252	&	37	&&  7040  &869	&	3928	&	244	\\
$P_{12}$ &	2	&&	58	&	42	&	18	&&	27	&	27	&	13	&&	23	&	3	&&	112	&	73	&	40	&	18	&&  3428  &385	&	1434	&	119	\\
$P_{12}$ &	5	&&	304	&	256	&	40	&&	71	&	70	&	37	&&	48	&	3	&&	426	&	262	&	164	&	35	&& 10736  &1688	&	4510	&	301	\\
$P_{12}$ &	10	&&	427	&	348	&	59	&&	132	&	131	&	76	&&	80	&	3	&&1079	&	627	&	451	&	64	&& 16267  &4092	&	7323	&	546	\\
$P_{12}$ &	50	&&3605	&3382	&	204	&&	241	&	240	&	157	&&	150	&	3	&&2611	&1381	&1230   &	76	&& 17738  &3904	&	8927	&	383	\\
$P_{12}$ &	100	&&2255	&2092	&	155	&&	268	&	269	&	184	&&	232	&	3	&&3861	&2006	&1854   &	77	&& 19525  &5250	&	9900	&	346	\\
$P_{12}$ &	200	&&2438	&2275	&	161	&&	497	&	495	&	366	&&	350	&	3	&&3758	&1964	&1793   &	87	&& 17479  &3969	&	8868	&	333	\\
$P_{13}$ &	2	&&	243	&	208	&	30	&&	32	&	32	&	17	&&	24	&	3	&&	223	&	138	&	85	&	25	&&  5209  &764	&	2190	&	174	\\
$P_{13}$ &	5	&&	755	&	676	&	63	&&	101	&	102	&	63	&&	73	&	4	&&	451	&	283	&	167	&	31	&& 16732  &1926	&	7233	&	335	\\
$P_{13}$ &	10	&&1465	&1357	&	89	&&	204	&	205	&	148	&&	92	&	4	&&	703	&	405	&	298	&	42	&& 26730  &3811	&	12255	&	714	\\
$P_{13}$ &	50	&&2401	&2277	&	108	&&	443	&	444	&	358	&&	182	&	4	&&1344	&	712	&	632	&	35  && 25503  &4475	&	12963	&	521	\\
$P_{13}$ &	100	&&3120	&2995	&	123	&&	389	&	389	&	301	&&	261	&	4	&&2257	&1172	&1086   &	39	&& 29193  &6021	&	14993	&	512	\\
$P_{13}$ &	200	&&5080	&4953	&	129	&&	485	&	486	&	351	&&	389	&	4	&&2247	&1168	&1078   &	41	&& 27111  &5128	&	13753	&	487	\\
$P_{14}$ &	2	&&	228	&	205	&	27	&&	27	&	19	&	16	&&	31	&	2	&&	297	&	185	&	112	&	30	&&  1244  &469	&	553	&	103	\\
$P_{14}$ &	5	&&	550	&	499	&	51	&&	49	&	38	&	26	&&	63	&	2	&&	654	&	419	&	235	&	43	&&   987  &533	&	513	&	128	\\
$P_{14}$ &	10	&&1182	&1065	&	98	&&	107	&	95	&	55	&&	109	&	2	&&1075	&	637	&	438	&	51  &&  1613  &843	&	953	&	213	\\
$P_{14}$ &	50	&&9765	&9509	&	246	&&	714	&	701	&	333	&&2791	&	2	&&4440	&2334	&2106   &	70	&& 23681  &10507	&	12094	&	416	\\
$P_{14}$ &	100	&&46252	&45820	&	436	&&3159	&3143	&1674	&&8770	&	2	&&6843	&3520	&3323   &	118	&& 27344  &12028	&	14018	&	511	\\
$P_{14}$ &	200	&&69437	&69109	&	325	&&3869	&3854	&2834	&&71000	&	4	&&6912	&3573	&3339   &	132	&& 18789  &7686	&	9767	&	426	\\
\hline
\end{tabular}
}
\end{center}
\label{NFNGDCCOMG2}
\end{table}

\begin{table}[h!]
\caption{Number of function and subgradient evaluations of DC components (Group 3).}
\begin{center}
\setlength\tabcolsep{6pt}
\resizebox{0.99\textwidth}{!}{
\begin{tabular}{@{}lrrrrrlrrrrrlrrrrrlrrrrrlrrrrrlrrrrr@{}}
\toprule
Prob.    & $n$&& \multicolumn{3}{c}{AggSub}&&\multicolumn{3}{c}{DBDC} &&\multicolumn{2}{c}{DCA} &&\multicolumn{4}{c}{ASM-DC} &&\multicolumn{4}{c}{TESGO}\\
\hline
 & &&$n_f$ &$n_{g_1}$ &$n_{g_2}$ &&$n_f$ &$n_{g_1}$ &$n_{g_2}$&& $n_{fg_1}$ &$n_{fg_2}$ &&$n_{f_1}$ &$n_{f_2}$ &$n_{g_1}$ &$n_{g_2}$ &&$n_{f_1}$&$n_{f_2}$ &$n_{g_1}$ &$n_{g_2}$    \\
 \cline{4-6} \cline{8-10} \cline{12-13} \cline{15-18} \cline{20-23}
$P_{15}$ &	2	&&	141	&	103	&	29	&&	70	&	65	&	43	&&	101	&	7	&&	290	&	191	&	99	&	36	&&	3016	&	 1148	&  1180	&	237	\\
$P_{16}$ &	2	&&	111	&	82	&	25	&&	23	&	21	&	16	&&	5	&	2	&&	215	&	135	&	80	&	28	&&	1642	&	 774	&   690	&	178	\\
$P_{16}$ &	5	&&	118	&	83	&	27	&&	22	&	22	&	16	&&	5	&	2	&&	262	&	158	&	104	&	30	&&	2254	&	 1007	&  1044	&	252	\\
$P_{16}$ &	10	&&	119	&	82	&	28	&&	18	&	18	&	14	&&	5	&	2	&&	323	&	194	&	130	&	32	&&	2926	&	 1312	&  1429	&	337	\\
$P_{16}$ &	50	&&	130	&	87	&	30	&&	26	&	25	&	20	&&	6	&	2	&&	539	&	301	&	238	&	33	&&	3250	&	 1484	&  2081	&	278	\\
$P_{16}$ &	100	&&	139	&	99	&	32	&&	14	&	15	&	11	&&	6	&	2	&&	846	&	451	&	395	&	32	&&	4227	&	 1889	&  2586	&	256	\\
$P_{16}$ &	200	&&	145	&	109	&	31	&&	18	&	19	&	19	&&	6	&	2	&&	955	&	513	&	442	&	35	&&	3653	&	 1543	&  2317	&	228	\\
$P_{17}$ &	2	&&	240	&	153	&	49	&&	62	&	50	&	34	&&	128	&	9	&&	352	&	235	&	116	&	36	&&	2123	&	 758	&   845	&	137	\\
$P_{18}$ &	2	&&	125	&	95	&	25	&&	48	&	42	&	28	&&	31	&	2	&&	274	&	179	&	94	&	34	&&	1003	&	 587	&   415	&	133	\\
$P_{18}$ &	5	&&	172	&	127	&	34	&&	147	&	140	&	79	&&	46	&	3	&&	452	&	286	&	166	&	44	&&	1386	&	 768	&   662	&	180	\\
$P_{18}$ &	10	&&	182	&	135	&	35	&&	165	&	160	&	89	&&	36	&	3	&&	533	&	326	&	207	&	46	&&	2169	&	 1149	&  1130	&	285	\\
$P_{18}$ &	50	&&	209	&	140	&	39	&&	180	&	181	&	95	&&	30	&	3	&&	780	&	453	&	327	&	46	&&	4334	&	 1843	&  2538	&	302	\\
$P_{18}$ &	100	&&	221	&	159	&	41	&&	161	&	162	&	88	&&	29	&	3	&&1038	&	582	&	456	&	47	&&	4320    &	1786	&  2584	&	263	\\
$P_{18}$ &	200	&&	239	&	181	&	43	&&	135	&	136	&	79	&&	27	&	3	&&1101	&	627	&	473	&	48	&&	3253    &	1197	&  2106	&	213	\\
$P_{19}$ &	2	&&	122	&	91	&	26	&&	56	&	53	&	32	&&	5	&	2	&&	247	&	158	&	88	&	31	&&	1322	&	 584	&   575	&	140	\\
$P_{20}$ &	2	&&	151	&	111	&	26	&&	23	&	14	&	13	&&	33	&	3	&&	258	&	165	&	93	&	28	&&	2213	&	 611	&   936	&	135	\\
$P_{20}$ &	5	&&	334	&	244	&	60	&&	448	&	429	&	287	&& 1087	&	26	&&	768	&	479	&	289	&	54	&&	15946	&	 3311	&  6645	&	419	\\
$P_{20}$ &	10	&& 3530	&  2397	&  758	&&15720	& 15660	& 11668	&& 1404	&	23	&& 1432	&  826	&	606	&	71	&&	27681	&	 6466	& 12363	&	662	\\
$P_{20}$ &	50	&&57539	& 51804	& 5550	&& 4796	&  4780	&  3553	&&29753	&  114	&&16395	& 8377	&8018	&	253	&&	70158	&	29780	& 34564	&	994	\\
$P_{20}$ &	100	&&68127	&64597	& 3491	&&33580	&33563	&20919	&&108401&	229	&&28946	&14629	&14317  &	281	&&	95306	&	42935	& 47362	&	927	\\
$P_{20}$ &	200	&&107855&105092	& 2760	&&12408	&12391	& 9157	&&396454&	89	&&19377	&9828	&9549	&	239	&&	60257	&	25700	& 29943	&	736	\\
\hline
\end{tabular}
}
\end{center}
\label{NFNGDCCOMG3}
\end{table}

Results presented in Tables \ref{NFNGDCCOMG1} --  \ref{NFNGDCCOMG3} show that with very few exceptions TESGO requires significantly more function and subgradient evaluations of both DC components than the other four local methods. As mentioned before, this is due to the fact that TESGO applies the local search method multiple times. However, taking into account that TESGO obtains higher quality solutions than other methods, the computational effort required by this method is reasonable.

\subsection{Comparison with global optimization solvers}
In this subsection, we present results for global minimization of the test problems from all three groups. The results of the proposed method is also compared with those obtained using well-known global optimization solvers BARON \cite{BARON} and LINDOGlobal \cite{LindoGlobal}. In addition, we consider \textit{two} hours time limit for solving each test problem. The results are given in Tables \ref{tableGROUP1}$-$\ref{tableGROUP3}, where we report the optimal value $f_{opt}$ obtained by a solver and errors $E_T, E_B$ and $E_L$ computed using \eqref{accuracy}. The notation $t_{lim}$ in the tables indicates that a solver reach the \textit{two} hours time limit. Since BARON solver is not applicable to all test problems, we use ``--" for such problems in tables.  This is due to the fact that BARON cannot handle the maximum function. If the maximum function is used in the first DC component $f_1$, then it can be rewritten without the maximum by introducing constraints and one extra variable. This is not possible if the second DC component $f_2$ has a maximum function, and thus these types of problems cannot be solved with BARON.

Here we apply the full version of the TESGO method. This means that we compute significantly more $\varepsilon$-subgradients of DC components in comparison with the simple version of TESGO. More specifically, we compute $m_1=\min\{100, 2n\}$ $\varepsilon$-subgradients of the first DC component and $m_2=\min\{30, 2n\}$ $\varepsilon$-subgradients of the second DC component. For the test problems with a large number of local minimizers, we set $m_1=\min\{150, 2n\}$ and $m_2=\min\{30, 2n\}$ (in tables these problems are indicated by $^*$). Finally, for some very complex problems, we set $m_1=\min\{200, 2n\}$ and $m_2=\min\{30, 2n\}$ (in tables these problems are indicated by $^{**}$).

\begin{table}[h!]
\caption{Results for TESGO, BARON, and LINDO (Group 1).}
\begin{center}
\setlength\tabcolsep{6pt}
\resizebox{0.99\textwidth}{!}{
\begin{tabular}{@{}lrrrrrrrrrrrrr@{}}
\toprule
Prob.       & $n$&& \multicolumn{3}{c}{TESGO}&&\multicolumn{3}{c}{BARON} &&\multicolumn{3}{c}{LINDO}\\
\hline
            &    && $f_{opt}$ & $E_T$ & CPU  &&  $f_{opt}$ & $E_B$ & CPU &&$f_{opt}$ & $E_L$ & CPU  \\
                    \cline{4-6} \cline{8-10} \cline{12-14}
$P_1$ &	2	&&	0.0000	&	0.0000	&	0.02	&&	0.0000	&	0.0000	&0.05	    &&	0.0000	&	0.0000	&	0.02	\\
$P_2$ &	4	&&	0.0000	&	0.0000	&	0.00	&&	0.0000	&	0.0000	&$t_{lim}$	&&	0.0000	&	0.0000	&	0.08	 \\
$P_3$ &	2	&&	0.5000	&	0.0000	&	0.00	&&	0.5000	&	0.0000	&0.06   	&&	0.5000	&	0.0000	&	0.13	\\
$P_4$ &	3	&&	3.5000	&	0.0000	&	0.00	&&	3.5000	&	0.0000	&0.22   	&&	3.5000	&	0.0000	&	0.08	\\
$P_5$ &	2	&&	-0.5000	&	0.0000	&	0.00	&&	-0.5000	&	0.0000	&0.11   	&&	-0.5000	&	0.0000	&	0.03	\\
$P_5$ &	5	&&	-3.5000	&	0.0000	&	0.20	&&	-3.5000	&	0.0000	&4.17   	&&	-3.5000	&	0.0000	&	0.09	\\
$P_5$ &	10	&&	-8.5000	&	0.0000	&	3.17	&&	-8.5000	&	0.0000	&$t_{lim}$	&&	-8.5000	&	0.0000	&   2.31    \\
$P_5$ &	50	&&-48.5000	&	0.0000	&	0.70	&&-47.5000	&	0.0202	&$t_{lim}$	&&-48.5000	&	0.0000	& $t_{lim}$	\\
$P_5$ &	100	&&-98.5000	&	0.0000	&	0.77	&&-90.5000	&	0.0804	&$t_{lim}$	&&-98.5000	&	0.0000	& $t_{lim}$	\\
$P_5$ &	200	&&-198.5000	&	0.0000	&	0.95	&&-182.5000	&	0.0802	&$t_{lim}$	&&-198.5000	&	0.0000	& $t_{lim}$	\\
$P_6$ &	3	&&116.3333	&	0.0000	&	0.00	&&116.3333	&	0.0000	&0.14   	&&116.3333	&	0.0000	&	0.03	\\
$P_7$ &	2	&&	0.0000	&	0.0000	&	0.00	&&	0.0000	&	0.0000	&0.04   	&&	0.0000	&	0.0000	&	0.03	\\
$P_7$ &	5	&&	0.0000	&	0.0000	&	0.00	&&	0.0000	&	0.0000	&0.04   	&&	0.0000	&	0.0000	&	0.22	\\
$P_7$ &	10	&&	0.0000	&	0.0000	&	0.02	&&	0.0000	&	0.0000	&0.04   	&&	0.0014	&	0.0014	&	1.93	\\
$P_7$ &	50	&&	0.0000	&	0.0000	&	0.34	&&	0.0000	&	0.0000	&$t_{lim}$	&&	0.0306	&	0.0306	&	 $t_{lim}$	\\
$P_7$ &	100	&&	0.0000	&	0.0000	&	0.86	&&	0.0000	&	0.0000	&0.39   	&&	0.0273	&	0.0273	&	$t_{lim}$	 \\
$P_7$ &	200	&&	0.0000	&	0.0000	&	3.97	&&	0.0000	&	0.0000	&2.43   	&&	0.0348	&	0.0348	&	$t_{lim}$	 \\
$P_8$ &	5	&&	0.0000	&	0.0000	&	0.08	&&  	--	&   	--	&  	--   	&&	0.0000	&	0.0000	&	263.74	\\
$P_8$ &	10	&&	0.0000	&	0.0000	&	0.44	&&  	--	&   	--	&  	--   	&&	0.0000	&	0.0000	&	$t_{lim}$	 \\
$P_8$ &	50	&&	0.0000	&	0.0000	&	72.17	&&  	--	&   	--	&  	--  	&&	0.0000	&	0.0000	&	4.40	\\
$P_8$ &	100	&&	0.0001	&	0.0001	&	1997.05	&&  	--	&   	--	&  	--  	&&	0.0000	&	0.0000	&	3.89	\\
$P_8$ &	200	&&	6.0400$^{**}$	&	6.0400$^{**}$	&	292.70$^{**}$	&& 	--	&   	--	&   	--	&&	0.0000	&	 0.0000	&	5.46	\\
\hline
\multicolumn{14}{l}{$^{**}$ in TESGO we have set $m_1=\min\{200, 2n\}$ and $m_2=\min\{30, 2n\}$}
\end{tabular}
}
\end{center}
\label{tableGROUP1}
\end{table}

\begin{table}[h!]
\caption{Results for TESGO, BARON, and LINDO (Group 2).}
\begin{center}
\setlength\tabcolsep{6pt}
\resizebox{0.99\textwidth}{!}{
\begin{tabular}{@{}lrrrrrrrrrrrrr@{}}
\toprule
Prob.       & $n$&& \multicolumn{3}{c}{TESGO}&&\multicolumn{3}{c}{BARON} &&\multicolumn{3}{c}{LINDO}\\
\hline
            &    && $f_{opt}$ & $E_T$ & CPU  &&  $f_{opt}$ & $E_B$ & CPU &&$f_{opt}$ & $E_L$ & CPU  \\
                    \cline{4-6} \cline{8-10} \cline{12-14}
$P_9$    &	2	&&-153.3333	&	0.0000	&	0.02	&&	--	&	--	&	--	&&	-153.3333	&	0.0000	&	376.69	\\
$P_9$  	 &	5	&&-436.6667	&	0.0000	&	0.03	&&	--	&	--	&	--	&&	-436.6667	&	0.0000	&	819.50	\\
$P_9$  	 &	10	&&-929.0906	&	0.0000	&	0.28	&&	--	&	 -- 	&	 -- 	&&	-929.0909	&	0.0000	&	565.11	 \\
$P_9$  	 &	50	&&-4921.9601&	0.0000	&	6.00	&&  -- 	&	 -- 	&	 -- 	&&	-4921.7032	&	0.0001	&	 $t_{lim}$	\\
$P_9$  	 &	100	&&-9920.9888$^{**}$	&	0.0000$^{**}$	&	2219.03$^{**}$	&&  -- 	&	 -- 	&	 -- 	&&	-9918.9905	 &	0.0002	&	$t_{lim}$	\\
$P_{10}$ &	2	&&-247.8125	&	0.0000	&	0.00	&&	-247.8125	&	0.0000	&	0.20	&&	-247.8125	&	0.0000	&	 137.26	\\
$P_{10}$ &	5	&&-578.4626	&	0.0000	&	0.02	&&	-578.4626	&	0.0000	&	0.31	&&	-578.4626	&	0.0000	&	 198.85	\\
$P_{10}$ &	10	&&-1006.8613$^*$	&	0.0000$^*$	&	0.17$^*$ 	&&	-1006.8616	&	0.0000	&	0.34	&&	-1006.8616	&	 0.0000	&	$t_{lim}$	\\
$P_{10}$ &	50	&&-3564.2274$^{**}$	&	0.0000$^{**}$	&163.52$^{**}$	&&	-3564.2275	&	0.0000	&	1.40	 &&-3538.2191&0.0073	&	$t_{lim}$	\\
$P_{10}$ &	100	&&-7297.9529&	0.0000	&	263.20	&&	-7297.9530	&	0.0000	&	5.99	&&	-7235.4972	&	0.0086	 &	$t_{lim}$	\\
$P_{11}$ &	2	&&	-0.5000	&	0.0000	&	0.00	&&	-0.5000	&	0.0000	&	0.12	&&	-0.5000	&	0.0000	&	0.04	 \\
$P_{11}$ &	5	&&	-3.5000	&	0.0000	&	0.16	&&	-3.5000	&	0.0000	&	4.17	&&	-3.5000	&	0.0000	&	0.10	 \\
$P_{11}$ &	10	&&	-8.5000	&	0.0000	&	0.16	&&	-8.5000	&	0.0000	&	5655.25	&&	-8.5000	&	0.0000	&	2.47	 \\
$P_{11}$ &	50	&&-48.5000	&	0.0000	&	0.70	&&	-47.5000	&	0.0202	&	$t_{lim}$	&&	-48.5000	&	 0.0000	&	$t_{lim}$	\\
$P_{11}$ &	100	&&-98.5000	&	0.0000	&	0.78	&&	-89.5000	&	0.0905	&	$t_{lim}$	&&	-98.5000	&	 0.0000	&	$t_{lim}$	\\
$P_{11}$ &	200	&&-198.5000	&	0.0000	&	0.95	&&	-180.5000	&	0.0902	&	$t_{lim}$	&&	-198.5000	&	 0.0000	&	$t_{lim}$	\\
$P_{12}$ &	2	&&	0.0000	&	0.0000	&	0.00	&&	 -- 	&	 -- 	&	 -- 	&&	0.0000	&	0.0000	&	0.21	 \\
$P_{12}$ &	5	&&	-1.8541	&	0.0000	&	0.05	&&	 -- 	&	 -- 	&	 -- 	&&	-1.8541	&	0.0000	&	22.18	 \\
$P_{12}$ &	10	&&	-4.9443	&	0.0000	&	0.97	&&	 -- 	&	 -- 	&	 -- 	&&	-4.9443	&	0.0000	&	 $t_{lim}$	\\
$P_{12}$ &	50	&&	-29.6656$^*$	&	0.0000$^*$	&	3.41$^*$	&&	 -- 	&	 -- 	&	 -- 	&&	-29.6656	&	 0.0000	&	$t_{lim}$	\\
$P_{12}$ &	100	&&	-60.5673$^{**}$	&	0.0000$^{**}$	&	579.05$^{**}$	&&	 -- 	&	 -- 	&	 -- 	&&	-60.5673	 &	0.0000	&	$t_{lim}$	\\
$P_{12}$ &	200	&&	-122.3707$^*$	&	0.0000$^*$	&	42.17$^*$	&&	 -- 	&	 -- 	&	 -- 	&&	-122.3707	&	 0.0000	&	$t_{lim}$	\\
$P_{13}$ &	2	&&	-5.0000	&	0.0000	&	0.61	&&	 -- 	&	 -- 	&	 -- 	&&	-5.0000	&	0.0000	&	0.14	 \\
$P_{13}$ &	5	&&-21.8541	&	0.0000	&	0.05	&&	 -- 	&	 -- 	&	 -- 	&&	-21.8541	&	0.0000	&	 1.66	\\
$P_{13}$ &	10	&&-48.9443	&	0.0196	&	0.08	&&	 -- 	&	 -- 	&	 -- 	&&	-49.9443	&	0.0000	&	 $t_{lim}$	\\
$P_{13}$ &	50	&&-226.0525$^{**}$&	0.1733$^{**}$	&	31.13$^{**}$	&&	 -- 	&	 -- 	&	 -- 	&&	 -273.6652	&	0.0000	&	$t_{lim}$	\\
$P_{13}$ &	100	&&-555.5672$^{**}$&	0.0000$^{**}$	&	560.44$^{**}$	&&	 -- 	&	 -- 	&	 -- 	&&	 -554.5672	&	0.0018	&	$t_{lim}$	\\
$P_{13}$ &	200	&&-939.2915$^{**}$&	0.1584$^{**}$	&	195.28$^{**}$	&&	 -- 	&	 -- 	&	 -- 	&&	 -1116.3273	&	0.0000	&	$t_{lim}$	\\
$P_{14}$ &	2	&&	-1.0000	&	0.0000	&	0.00	&&	 -- 	&	 -- 	&	 -- 	&&	-1.0000	&	0.0000	&	0.34	 \\
$P_{14}$ &	5	&&	-3.4167	&	0.0000	&	0.00	&&	 -- 	&	 -- 	&	 -- 	&&	-3.4167	&	0.0000	&	15.60	 \\
$P_{14}$ &	10	&&-11.2897	&	0.0000	&	0.03	&&	 -- 	&	 -- 	&	 -- 	&&	-11.2897	&	0.0000	&	 127.92	\\
$P_{14}$ &	50	&&-126.9603	&	0.0000	&	6.50	&&	 -- 	&	 -- 	&	 -- 	&&	-125.8108	&	0.0090	&	 $t_{lim}$	\\
$P_{14}$ &	100	&&-320.7221	&	0.0000	&	31.81	&&	 -- 	&	 -- 	&	 -- 	&&	-311.7939	&	0.0278	&	 $t_{lim}$	\\
$P_{14}$ &	200	&&-776.3854	&	0.0016	&	25.84	&&	 -- 	&	 -- 	&	 -- 	&&	-774.7685	&	0.0036	&	 $t_{lim}$	\\
\hline
\multicolumn{14}{l}{$^{*}$ in TESGO we have set $m_1=\min\{150, 2n\}$ and $m_2=\min\{30, 2n\}$}\\
\multicolumn{14}{l}{$^{**}$ in TESGO we have set $m_1=\min\{200, 2n\}$ and $m_2=\min\{30, 2n\}$}
\end{tabular}
}
\end{center}
\label{tableGROUP2}
\end{table}

\begin{table}[h!]
\caption{Results for TESGO, BARON, and LINDO (Group 3).}
\begin{center}
\setlength\tabcolsep{6pt}
\resizebox{0.99\textwidth}{!}{
\begin{tabular}{@{}lrrrrrrrrrrrrrrrrr@{}}
\toprule
Prob.       & $n$&& \multicolumn{3}{c}{TESGO}&&\multicolumn{3}{c}{BARON} &&\multicolumn{3}{c}{LINDO}\\
\hline
            &    && $f_{opt}$ & $E_T$ & CPU  &&  $f_{opt}$ & $E_B$ & CPU &&$f_{opt}$ & $E_L$ & CPU  \\
                    \cline{4-6} \cline{8-10} \cline{12-14}
$P_{15}$	&	2	&&	-0.3524	&	0.0000	&	0.00	&&	-0.3524	&	0.0000	&	0.06	&&	-0.3524	&	0.0000	&	0.48	 \\
$P_{16}$	&	2	&&	0.0000	&	0.0000	&	0.02	&&	0.0000	&	0.0000	&	0.16	&&	0.0000	&	0.0000	&	0.04	 \\
$P_{16}$	&	5	&&	0.0000	&	0.0000	&	0.00	&&	0.0000	&	0.0000	&	2.73	&&	0.0000	&	0.0000	&	0.11	 \\
$P_{16}$	&	10	&&	0.0000	&	0.0000	&	0.02	&&	0.0000	&	0.0000	&	289.58	&&	0.0000	&	0.0000	&	0.16	 \\
$P_{16}$	&	50	&&	0.0000	&	0.0000	&	0.58	&&	0.0000	&	0.0000	&	$t_{lim}$	&&	0.0000	&	0.0000	&	 $t_{lim}$	\\
$P_{16}$	&	100	&&	0.0000	&	0.0000	&	0.94	&&	0.0000	&	0.0000	&	$t_{lim}$	&&	0.0000	&	0.0000	&	 $t_{lim}$	\\
$P_{16}$	&	200	&&	0.0000	&	0.0000	&	1.20	&&	0.0000	&	0.0000	&	$t_{lim}$	&&	0.0000	&	0.0000	&	 $t_{lim}$	\\
$P_{17}$	&	2	&&	-0.8333	&	0.0000	&	0.00	&&	-0.8333	&	0.0000	&	0.18	&&	-0.8333	&	0.0000	&	3.13	 \\
$P_{18}$	&	2	&&	-0.3750	&	0.0000	&	0.00	&&	-0.3750	&	0.0000	&	0.12	&&	-0.3750	&	0.0000	&	0.05	 \\
$P_{18}$	&	5	&&	-1.3750	&	0.0000	&	0.00	&&	-1.3750	&	0.0000	&	2.10	&&	-1.3750	&	0.0000	&	0.12	 \\
$P_{18}$	&	10	&&	-3.0417	&	0.0000	&	0.02	&&	-3.0417	&	0.0000	&	668.42	&&	-3.0417	&	0.0000	&	3.02	 \\
$P_{18}$	&	50	&&-16.3750	&	0.0000	&	0.75	&&-16.3750	&	0.0000	&	$t_{lim}$	&&-16.3750	&	0.0000	&	 $t_{lim}$	\\
$P_{18}$	&	100	&&-33.0417	&	0.0000	&	0.86	&&-33.0417	&	0.0000	&	$t_{lim}$	&&-33.0417	&	0.0000	&	 $t_{lim}$	\\
$P_{18}$	&	200	&&-66.3750	&	0.0000	&	0.95	&&-66.3750	&	0.0000	&	$t_{lim}$	&&-66.3750	&	0.0000	&	 $t_{lim}$	\\
$P_{19}$	&	2	&&	-0.2500	&	0.0000	&	0.00	&&	-0.2500	&	0.0000	&	0.11	&&	-0.2500	&	0.0000	&	0.04	 \\
$P_{20}$	&	2	&&	0.0000	&	0.0000	&	0.00	&&	0.0000	&	0.0000	&	0.04	&&	0.0000	&	0.0000	&	0.08	 \\
$P_{20}$	&	5	&&	0.0000	&	0.0000	&	0.06	&&	0.0000	&	0.0000	&	0.05	&&	0.0000	&	0.0000	&	0.06	 \\
$P_{20}$	&	10	&&	0.0000	&	0.0000	&	0.17	&&	0.0000	&	0.0000	&	0.04	&&	0.0000	&	0.0000	&	0.11	 \\
$P_{20}$	&	50	&&	0.0000	&	0.0000	&	7.88	&&	0.0000	&	0.0000	&	0.06	&&	0.0000	&	0.0000	&	2.08	 \\
$P_{20}$	&	100	&&	0.0000	&	0.0000	&	162.97	&&	0.0000	&	0.0000	&	0.21	&&	0.0000	&	0.0000	&	5.52	 \\
$P_{20}$	&	200	&&	0.0000	&	0.0000	&	76.59	&&	0.0000	&	0.0000	&	0.64	&&	0.0000	&	0.0000	&	16.40	 \\
\hline
\end{tabular}
}
\end{center}
\label{tableGROUP3}
\end{table}

Results from Tables \ref{tableGROUP1}$-$\ref{tableGROUP3} show that TESGO  finds global solutions in 72 cases out of 77, the BARON solver in 43 cases out of 49 and the LINDOGlobal solver in 71 cases out of 77. However, both BARON and LINDOGlobal require significantly more CPU time than TESGO. The only exceptions are $P_8$ and $P_{20}$ with $n=50, 100, 200$. In many cases, BARON and LINDOGlobal are forced to stop due to the two hours time limit. These results clearly indicate that, in general, TESGO is able to find accurate solutions to most DC optimization problems by using significantly less computational effort than BARON and LINDOGlobal.

Finally, in Table \ref{NfNg-GODC}, we report the number of function and subgradient evaluations required by  TESGO to solve DC optimization problems to global optimality. These numbers are computed as an average of the number of function and subgradient evaluations of DC components. We only report these results for TESGO as such information for BARON and LINDOGlobal cannot be extracted.  We can see from this table that in most cases the TESGO method uses a reasonable number of function and subgradient evaluations. Depending on the starting point, the large number of local minimizers can lead to large number of function and subgradient evaluations. Problems $P_8$ with $n=50, 100, 200$, $P_9$ with $n = 100$, $P_{10}$ with $n=50, 100$, $P_{12}$ with $n=100$ and $P_{13}$ with $n=50, 100, 200$ are among such problems. In these problems, TESGO requires a large number of function and subgradient evaluations.

\begin{table}[h!]
\caption{Number of function and subgradient evaluations for TESGO}
\begin{center}
\setlength\tabcolsep{6pt}
\resizebox{0.99\textwidth}{!}{
\begin{tabular}{@{}lrrrrlrrrlrrrr@{}}
\toprule
Prob. &$n$ & $n_f$  & $n_g$  &&  Prob. & $n$ & $n_f$  & $n_g$  && Prob. & $n$ & $n_f$  & $n_g$    \\
\hline
\multicolumn{14}{c}{Group 1 }\\
\hline
$P_1$	&	2	&	225	    &	159	    &&	$P_5$	&	100	&	3727	&	1997	&&	$P_7$	&	200	&	9880	&	4537	 \\
$P_2$	&	4	&	3631	&	1452	&&	$P_5$	&	200	&	4133	&	2094	&&	$P_8$	&	5	&	29000	&	10147	 \\
$P_3$	&	2	&	8428	&	2767	&&	$P_6$	&	3	&	266	    &	196	    &&	$P_8$	&	10	&	66402	&	23332	 \\
$P_4$	&	3	&	1983	&	771	    &&	$P_7$	&	2	&	1257	&	532	    &&	$P_8$	&	50	&	147519	&	50830	 \\
$P_5$	&	2	&	2381	&	855	    &&	$P_7$	&	5	&	5791	&	2043	&&	$P_8$	&	100	&	140322	&	49153	 \\
$P_5$	&	5	&	3239	&	1389	&&	$P_7$	&	10	&	9917	&	3728	&&	$P_8$	&	200	&	102488$^{**}$	&	 39073$^{**}$ \\
$P_5$	&	10	&	4806	&	1991	&&	$P_7$	&	50	&	9655	&	4442	&&								\\
$P_5$	&	50	&	3420	&	1883	&&	$P_7$	&	100	&	7243	&	3521	&&								\\
\hline
\multicolumn{14}{c}{Group 2}\\
\hline																							
$P_9$	&	2	&	3701	&	2432	&&	$P_{11}$	&	10	&	4646	&	3632    &&	$P_{13}$	&	10	&	14626	&	 9409	\\
$P_9$	&	5	&	13261	&	9196	&&	$P_{11}$	&	50	&	2593	&	2240    &&	$P_{13}$	&	50	&	 107211$^{**}$	&	93647$^{**}$ \\
$P_9$	&	10	&	12234	&	7364	&&	$P_{11}$	&	100	&	2968	&	2656    &&	$P_{13}$	&	100	&	 135975$^{**}$	&	115026$^{**}$ \\
$P_9$	&	50	&	13273	&	12370	&&	$P_{11}$	&	200	&	3183	&	2981	&&	$P_{13}$	&	200	&	69674$^{**}$	 &	55864$^{**}$ \\
$P_9$	&	100	&325787$^{**}$&	324040$^{**}$ && $P_{12}$	&	2	&	2260	&	1443&&	$P_{14}$	&	 2	&	337	&	364	\\
$P_{10}$&	2	&	3444	&	2319    &&	$P_{12}$	&	5	&	11406	&	7401	        &&	 $P_{14}$	&	5	&	353	&	478	\\
$P_{10}$&	5	&	8177	&	5808	                     &&	$P_{12}$	&	10	&	17760	&	12918	        &&	 $P_{14}$	&	10	&	511	&	761	\\
$P_{10}$&	10	&25962$^*$  &	17849$^*$	             &&	$P_{12}$	&	50	&	28432$^{*}$	&	17843$^{*}$	&&	 $P_{14}$	&	50	&	13300	&	12397	\\
$P_{10}$&	50	&330857$^{**}$&329517$^{**}$&&	$P_{12}$	&	100	&	110911$^{**}$	&	96628$^{**}$	 &&	$P_{14}$	&	100	&	13216	&	11688	\\
$P_{10}$&	100	&	62905	&	61873	&&	$P_{12}$	&	200	&	32535$^{*}$	&	20883$^{*}$	&&	$P_{14}$	&	200	&	 10735	&	9685	\\
$P_{11}$&	2	&	1969	&	1589	&&	$P_{13}$	&	2	&	4589	&	2952	&&								\\
$P_{11}$&	5	&	2916	&	2026	&&	$P_{13}$	&	5	&	13040	&	8531	&&								\\
\hline
\multicolumn{14}{c}{Group 3}\\
\hline	
$P_{15}$	&	2	&	5151	&	1780	&&	$P_{17}$	&	2	&	1209	&	428	    &&	$P_{19}$	&	2	&	719	    &	304	\\
$P_{16}$	&	2	&	2172	&	787	    &&	$P_{18}$	&	2	&	408	    &	194	    &&	$P_{20}$	&	2	&	565   	 &	282	\\
$P_{16}$	&	5	&	981	    &	518	    &&	$P_{18}$	&	5	&	1515	&	627	    &&	$P_{20}$	&	5	&	15638	 &	6023	\\
$P_{16}$	&	10	&	3802	&	1606	&&	$P_{18}$	&	10	&	4791	&	1876	&&	$P_{20}$	&	10	&	27899	 &	10966	\\
$P_{16}$	&	50	&	3122	&	1733	&&	$P_{18}$	&	50	&	3357	&	1901	&&	$P_{20}$	&	50	&	41399	 &	16070	\\
$P_{16}$	&	100	&	2951	&	1804	&&	$P_{18}$	&	100	&	3101	&	1827	&&	$P_{20}$	&	100	&	76660	 &	27748	\\
$P_{16}$	&	200	&	4009	&	2052	&&	$P_{18}$	&	200	&	2350	&	1593	&&	$P_{20}$	&	200	&	49536	 &	18500	\\
\hline
\multicolumn{14}{l}{$^{*}$ in TESGO we have set $m_1=\min\{150, 2n\}$ and $m_2=\min\{30, 2n\}$}\\
\multicolumn{14}{l}{$^{**}$ in TESGO we have set $m_1=\min\{200, 2n\}$ and $m_2=\min\{30, 2n\}$}
\end{tabular}
}
\end{center}
\label{NfNg-GODC}
\end{table}

\section{Conclusions} \label{conclusions}
In this paper, a new algorithm, the truncated $\varepsilon$-subdifferential method, is developed to globally minimize DC functions subject to box-constraints. It is a hybrid method based on the combination of a local search and a special procedure for escaping from solutions of a DC function which are not global minimizers. A local search method is applied to find a stationary point (in our case a critical point) of the DC optimization problem. Then the escaping procedure is employed to escape from this point by finding a better initial point for a local search.

We compute subsets of the $\varepsilon$-subdifferentials of DC components.
Then we calculate the deviation of the subset of the $\varepsilon$-subdifferential of the second DC component from the subset of the $\varepsilon$-subdifferential of the first DC component. If this deviation is positive then we utilize the $\varepsilon$-subgradient of the second DC component providing this deviation  to formulate a subproblem with a convex objective function. The solution to this subproblem is used as a starting point for a local search. The convergence of the conceptual version of the proposed method is studied and its implementation is discussed in detail.

The performance of the new method is demonstrated using a large number of academic test problems for DC optimization. Based on extensive numerical results it is shown that the proposed method is able to significantly improve the quality of solutions obtained by a local method using limited computational effort. In addition, we apply the developed method to find global solutions to DC optimization problems. Results show that the new method is able to find global solutions by increasing the number of $\varepsilon$-subgradients calculations in the escaping procedure. Comparison with two widely used global optimization solvers shows that the proposed method is efficient and accurate for solving DC optimization problems to global optimality using significantly less computational effort.

\section*{Statements and Declarations}

\bibliographystyle{spmpsci}      

\bibliography{GOepsilon}   

\begin{thebibliography}{10}
\providecommand{\url}[1]{{#1}}
\providecommand{\urlprefix}{URL }
\expandafter\ifx\csname urlstyle\endcsname\relax
  \providecommand{\doi}[1]{DOI~\discretionary{}{}{}#1}\else
  \providecommand{\doi}{DOI~\discretionary{}{}{}\begingroup
  \urlstyle{rm}\Url}\fi

\bibitem{Ackooij2021}
Ackooij, W., Demassey, S., Javal, P., Morais, H., de~Oliveira, W., Swaminathan,
  B.: A bundle method for nonsmooth {DC} programming with application to
  chance-constrained problems.
\newblock Comput. Optim. Appl. \textbf{78}(1), 451--490 (2021).
\newblock \doi{https://doi.org/10.1007/s10589-020-00241-8}

\bibitem{Antao2005}
An, L.T.H., Tao, P.D.: The {DC} (difference of convex functions) programming
  and {DCA} revisited with {DC} models of real world nonconvex optimization
  problems.
\newblock Ann. Oper. Res. \textbf{133}, 23--46 (2005).
\newblock \doi{https://doi.org/10.1007/s10479-004-5022-1}

\bibitem{Antao2012}
An, L.T.H., Tao, P.D., Ngai, H.V.: Exact penalty and error bounds in {DC}
  programming.
\newblock J. Glob. Optim. \textbf{52}(3), 509--535 (2012).
\newblock \doi{https://doi.org/10.1007/s10898-011-9765-3}

\bibitem{Artacho2020}
Artacho, F.J.A., Campoy, R., Vuong, P.T.: Using positive spanning sets to
  achieve $d$-stationarity with the boosted {DC} algorithm.
\newblock Vietnam J. Math. \textbf{48}(2), 363--376 (2020).
\newblock \doi{https://doi.org/10.1007/s10013-020-00400-8}

\bibitem{Artacho2018}
Artacho, F.J.A., Fleming, R.M.T., Vuong, P.T.: Accelerating the {DC} algorithm
  for smooth functions.
\newblock Math. Program. \textbf{169}, 95--118 (2018).
\newblock \doi{https://doi.org/10.1007/s10107-017-1180-1}

\bibitem{AV2020}
Artacho, F.J.A., Vuong, P.T.: The boosted difference of convex functions
  algorithm for nonsmooth functions.
\newblock SIAM J. Optim. \textbf{30}(1), 980--1006 (2020).
\newblock \doi{https://doi.org/10.1137/18M123339X}

\bibitem{Bagtahhos2021}
Bagirov, A.M., Hoseini~Monjezi, N., Taheri, S.: An augmented subgradient method
  for minimizing nonsmooth {DC} functions.
\newblock Comput. Optim. Appl. \textbf{80}(1), 411--438 (2021).
\newblock \doi{https://doi.org/10.1007/s10589-021-00304-4}

\bibitem{BagKarMak2014}
Bagirov, A.M., Karmitsa, N., M\"akel\"a, M.M.: Introduction to Nonsmooth
  Optimization.
\newblock Springer, Cham (2014).
\newblock \doi{https://doi.org/10.1007/978-3-319-08114-4}

\bibitem{BagKarTah2020}
Bagirov, A.M., Karmitsa, N., Taheri, S.: Partitional Clustering via Nonsmooth
  Optimization.
\newblock Springer, Cham (2020).
\newblock \doi{https://doi.org/10.1007/978-3-030-37826-4}

\bibitem{Bagirov2021}
Bagirov, A.M., Taheri, S., Cimen, E.: Incremental {DC} optimization algorithm
  for large-scale clusterwise linear regression.
\newblock J. Comput. Appl. Math \textbf{389}, 113323 (2021).
\newblock \doi{https://doi.org/10.1016/j.cam.2020.113323}

\bibitem{TUCreport2019}
Bagirov, A.M., Taheri, S., Joki, K., Karmitsa, N., M\"akel\"a, M.M.: A new
  subgradient based method for nonsmooth {DC} programming, {TUCS}.
\newblock Tech. rep., No. 1201, Turku Centre for Computer Science, Turku (2019)

\bibitem{Bagtahjokkarmak2021}
Bagirov, A.M., Taheri, S., Joki, K., Karmitsa, N., M\"akel\"a, M.M.: Aggregate
  subgradient method for nonsmooth {DC} optimization.
\newblock Optim. Lett. \textbf{15}(1), 83--96 (2021).
\newblock \doi{https://doi.org/10.1007/s11590-020-01586-z}

\bibitem{bagi2016}
Bagirov, A.M., Taheri, S., Ugon, J.: Nonsmooth {DC} programming approach to the
  minimum sum-of-squares clustering problems.
\newblock Pattern Recognit. \textbf{53}, 12--24 (2016).
\newblock \doi{https://doi.org/10.1016/j.patcog.2015.11.011}

\bibitem{Bagirov2011}
Bagirov, A.M., Ugon, J.: Codifferential method for minimizing nonsmooth {DC}
  functions.
\newblock J. Glob. Optim. \textbf{50}, 3--22 (2011).
\newblock \doi{https://doi.org/10.1007/s10898-010-9569-x}

\bibitem{Clarke1983}
Clarke, F.H.: Optimization and Nonsmooth Analysis.
\newblock Wiley-Interscience, New York (1983)

\bibitem{czyzyketal1998}
Czyzyk, J., Mesnier, M.P., Mor{\'{e}}, J.J.: The {NEOS} {S}erver.
\newblock IEEE Comput. Sci. Eng. \textbf{5}(3), 68--75 (1998).
\newblock \doi{https://doi.org/10.1109/99.714603}

\bibitem{Demyanov1982}
Demyanov, V.F., Vasilyev, L.: Nondifferentiable Optimization.
\newblock Optimization Software, New York (1986)

\bibitem{dolan2001}
Dolan, E.D.: The {NEOS} {S}erver 4.0 administrative guide.
\newblock Technical Memorandum ANL/MCS-TM-250, Mathematics and Computer Science
  Division, Argonne National Laboratory (2001)

\bibitem{Dolgopolik2018}
Dolgopolik, M.V.: A convergence analysis of the method of codifferential
  descent.
\newblock Comput. Optim. Appl. \textbf{71}(1), 879--913 (2018).
\newblock \doi{https://doi.org/10.1007/s10589-018-0024-0}

\bibitem{Ferrer2015}
Ferrer, A., Bagirov, A.M., Beliakov, G.: Solving {DC} programs using the
  cutting angle method.
\newblock J. Glob. Optim. \textbf{61}, 71--89 (2015).
\newblock \doi{https://doi.org/10.1007/s10898-014-0159-1}

\bibitem{Frangioni1996}
Frangioni, A.: Solving semidefinite quadratic problems within nonsmooth
  optimization algorithms.
\newblock Computers \& Oper. Res. \textbf{23}(11), 1099--1118 (1996).
\newblock \doi{https://doi.org/10.1016/0305-0548(96)00006-8}

\bibitem{Gaudioso2018}
Gaudioso, M., Giallombardo, G., Miglionico, G., Bagirov, A.M.: Minimizing
  nonsmooth {DC} functions via successive {DC} piecewise-affine approximations.
\newblock J. Glob. Optim. \textbf{71}(1), 37--55 (2018).
\newblock \doi{https://doi.org/10.1007/s10898-017-0568-z}

\bibitem{Goldstein1977}
Goldstein, A.A.: Optimization of lipschitz continuous functions.
\newblock Math. Program. \textbf{13}, 14--22 (1977).
\newblock \doi{https://doi.org/10.1007/BF01584320}

\bibitem{groppmore1997}
Gropp, W., Mor{\'{e}}, J.J.: Optimization environments and the {NEOS} {S}erver.
\newblock In: M.D. Buhman, A.~Iserles (eds.) Approximation Theory and
  Optimization, pp. 167--182. Cambridge University Press (1997)

\bibitem{Hiriart1988}
Hiriart-Urruty, J.B.: From convex optimization to nonconvex optimization.
  {N}ecessary and sufficient conditions for global optimality.
\newblock In: F.H. Clarke, V.F. Dem’yanov, F.~Giannessi (eds.) Nonsmooth
  Optimization and Related Topics, Ettore Majorana International Sciences
  Series 43, pp. 219--239. Springer, Boston (1989)

\bibitem{Horst1995}
Horst, R., Pardalos, P.M., Thoai, N.V.: Introduction to Global Optimization.
\newblock Kluwer Academic Publishers, Dordrecht (1995)

\bibitem{Horst1999}
Horst, R., Thoai, N.V.: {DC} programming: {O}verview.
\newblock J. Optim. Theory Appl. \textbf{103}(1), 1--43 (1999).
\newblock \doi{https://doi.org/10.1023/A:1021765131316}

\bibitem{Horst1991}
Horst, R., Thoai, N.V., Benson, H.P.: Concave minimization via conical
  partitions and polyhedral outer approximation.
\newblock Math. Program. \textbf{50}, 259--274 (1991).
\newblock \doi{https://doi.org/10.1007/BF01594938}

\bibitem{Horst1987}
Horst, R., Thoai, N.V., Tuy, H.: Outer approximation by polyhedral convex sets.
\newblock Oper.-Res.-Spektrum \textbf{9}, 153--159 (1987).
\newblock \doi{https://doi.org/10.1007/BF01721096}

\bibitem{Horst1996}
Horst, R., Tuy, H.: Global Optimization (Deterministic Approach).
\newblock Springer Verlag, Berlin, Germany (1996)

\bibitem{Joki2017}
Joki, K., Bagirov, A.M., Karmitsa, N., M\"akel\"a, M.M.: A proximal bundle
  method for nonsmooth {DC} optimization utilizing nonconvex cutting planes.
\newblock J. Glob. Optim. \textbf{68}(1), 501--535 (2017).
\newblock \doi{https://doi.org/10.1007/s10898-016-0488-3}

\bibitem{JokBagKarMakTah2018}
Joki, K., Bagirov, A.M., Karmitsa, N., M\"akel\"a, M.M., Taheri, S.: Double
  bundle method for finding {C}larke stationary points in nonsmooth {DC}
  programming.
\newblock SIAM J. Optim. \textbf{28}(2), 1892--1919 (2018).
\newblock \doi{https://doi.org/10.1137/16M1115733}

\bibitem{Kiwiel1989}
Kiwiel, K.C.: A dual method for certain positive semidefinite quadratic
  programming problems.
\newblock SIAM J. Sci. Stat. Comput. \textbf{10}(1), 175--186 (1989).
\newblock \doi{https://doi.org/10.1137/0910013}

\bibitem{LindoGlobal}
Lin, Y., Schrage, L.: The global solver in the {LINDO} {API}.
\newblock Optim. Methods Softw. \textbf{24}(4-5), 657--668 (2009).
\newblock \doi{https://doi.org/10.1080/10556780902753221}

\bibitem{MakNei1992}
M\"akel\"a, M.M., Neittaanm\"aki, P.: Nonsmooth Optimization: Analysis and
  Algorithms with Applications to Optimal Control.
\newblock World Scientific Publishing Co, Singapore (1992)

\bibitem{Dolan2002}
Mor{\'{e}}, J., Dolan, E.: Benchmarking optimization software with performance
  profiles.
\newblock Math. Program. \textbf{91}, 201--213 (2002).
\newblock \doi{https://doi.org/10.1007/s101070100263}

\bibitem{Nurminski2008}
Nurminski, E.A.: Projection onto polyhedra in outer representation.
\newblock Comput. Math. \& Math. Phys. \textbf{48}(3), 367--375 (2008).
\newblock \doi{https://doi.org/10.1134/S0965542508030044}

\bibitem{Oliveira2019}
de~Oliveira, W.: Proximal bundle methods for nonsmooth {DC} programming.
\newblock J. Glob. Optim. \textbf{75}, 523--563 (2019).
\newblock \doi{https://doi.org/10.1007/s10898-019-00755-4}

\bibitem{Oliveira2020}
de~Oliveira, W.: The {ABC} of {DC} programming.
\newblock Set-Valued Var. Anal. \textbf{28}(1), 679--706 (2020).
\newblock \doi{https://doi.org/10.1007/s11228-020-00566-w}

\bibitem{pinter1996}
Pinter, J.: Global Optimization in Action.
\newblock Kluwer Academic Publishers, Dordrecht (1996)

\bibitem{Rockafellar1970}
Rockafellar, R.T.: Convex Analysis.
\newblock Princeton University Press, Princeton (1970)

\bibitem{BARON}
Sahinidis, N.V.: {BARON} 23.5.23: {G}lobal Optimization of Mixed-Integer
  Nonlinear Programs, {\em User's Manual} (2023)

\bibitem{Tuy1998}
Tuy, H.: Convex Analysis and Global Optimization.
\newblock Kluwer Academic Publishers, Dordrecht (1998)

\bibitem{Wolfe1976}
Wolfe, P.H.: Finding the nearest point in a polytope.
\newblock Math. Program. \textbf{11}(2), 128--149 (1976).
\newblock \doi{https://doi.org/10.1007/BF01580381}

\end{thebibliography}

\section*{Appendix: Test problems}
All objective functions are DC functions:
$$
f(\x) =f_1(\x) - f_2(\x)
$$

\begin{problem}
\noindent \emph{DC version of Aluffi-Pentini's problem}
$$
f_1(\x) =0.25x_1^4 + 0.1x_1 + 0.5x_2^2,~~~~f_2(\x) = 0.5x_1^2
$$
$$
\x = (x_1, x_2) \in \R^2,~~x_i \in [-10, 10],~i=1,2.
$$
\end{problem}

\begin{problem}
\noindent \emph{Generalized DC Becker and Lago problem}
$$
f_1(\x) = \sum_{i=1}^n x_i^2 + 25n,~~~~f_2(\x) = 10 \sum_{i=1}^n |x_i|
$$
$$
\x \in \R^n,~~x_i \in [-10, 10],~i=1,\ldots,n.
$$
\end{problem}

\begin{problem}
\noindent \emph{Modified DC Camel Back problem}
$$
f_1(\x) = \frac{1}{6} + x_1^6 + 4x_1^2 + 4x_2^4 + |x_1|,~~~~f_2(\x) = 2.1x_1^4 + 4x_2^2
$$
$$
\x = (x_1, x_2) \in \R^2,~~x_i \in [-5, 5],~i=1,2.
$$
\end{problem}

\begin{problem}
\noindent \emph{}
$$
f_1(\x) = \sum_{i=2}^n \big((x_i - 1)^2 +x_{i-1}^2+x_i^2\big),~~~~f_2(\x) = \sum_{i=2}^n |x_{i-1} + x_i|
$$
$$
\x \in \R^n,~~x_i \in [-n, n],~i=1,\ldots,n.
$$
\end{problem}

\begin{problem}
\noindent \emph{}
$$
f_1(\x) = 2(x_1^2 + x_2^2),~~~~f_2(\x) = |x_1 + x_2|
$$
$$
\x \in \R^2,~~x_i \in [-10, 10],~i=1,2.
$$
\end{problem}

\begin{problem}
\noindent \emph{}
$$
f_1(\x) = 2\sum_{i=1}^{n-1} \max \Big\{x_{i+1} - x_i + 1, x_i^2 \Big\},~~~f_2(\x) = \sum_{i=1}^{n-1} \Big(x_i^2 +x_{i+1} - x_i +1\Big)
$$
$$
\x \in \R^n,~~x_i \in [-10, 10],~i=1,\ldots,n.
$$
\end{problem}

\end{document}